\documentstyle[amscd,amssymb,verbatim,12pt]{amsart}

\font\rm=cmr12

  \newcommand {\C}{{\mathbb C}}
\newcommand {\N}{{\mathbb N}}
  \newcommand {\R}{{\mathbb R}}
  \newcommand {\Z}{{\mathbb Z}}

  \newcommand {\B}{{\mathcal B}}
\renewcommand {\H}{{\mathcal H}}
  \newcommand {\M}{{\mathcal M}}

  \newcommand {\K}{{\mathcal K}}

\newcommand{\supp}{\operatorname{supp}}

\theoremstyle{plain}
\newtheorem{lem}{Lemma}[section]
\newtheorem{theo}[lem] {Theorem}
\newtheorem{cor}[lem]{Corollary}
\newtheorem{Def}[lem]{Definition}
\newtheorem{satz}[lem]{Proposition}

\newtheorem{bem}[lem]{Remark}

\theoremstyle{definition}

\theoremstyle{remark}

\begin{document}
\title[Index theory for coverings]
{\bf Index theory for coverings}
\date{June 2008}
\author{Boris Vaillant}
\address{bvaillant@@quantitative-consulting.eu}

\maketitle
\tableofcontents
\contentsline{section}
{\numberline{}}{}

\newpage

\section{Introduction}
\setcounter{equation}{0}

This paper is a translated and revised version of the author's Diploma thesis
\cite{Va}. The starting point is Atiyah's work \cite{A} on the index theory 
for elliptic operators on coverings of compact manifolds $M$.

Every elliptic differential operator $P\in \mbox{\rm Diff}^k(M;E,F)$ on sections of
hermitian vector bundles $E$, $F$ over a compact manifold is invertible upto a smoothing
operator and has a well-defined index. 





In the case of a regular covering $\overline M\to M$ with
an infinite group of deck transformations $\Gamma$ and lifted vector bundles $\overline E,\overline F,\overline P$,
the lifted operator $\overline P$ can have essential 
spectrum down to  $0$ and is in general not Fredholm. 
It is however still possible to define the $\Gamma$-{\it trace} of a 
$\Gamma$-equivariant operator $K$ on $\overline M$ as the integral of its local trace over a 
fundamental domain ${\mathcal F}$

$$\mbox{\rm tr}_\Gamma(K):=\int_{\mathcal F}\mbox{\rm tr}([K](x,x))dx.$$

The operator $\overline P$ is then  '$\Gamma$-Fredholm' and has a well defined
$\Gamma$-index

$$\mbox{\rm ind}_\Gamma(\overline P)=\mbox{\rm tr}_\Gamma(N(\overline P))-\mbox{\rm tr}_\Gamma(N(\overline{P^*})).$$

It is shown in \cite{A} that


\begin{equation}\label{4}
\mbox{ind}_\Gamma(\overline P)=\mbox{ind}(P).
\end{equation}

In this paper, we will look at $\Gamma$-index Theorems for coverings $\overline N\to N,$ 
where the base manifold $N$ is noncompact but still has bounded geometry. 





The simplest case of a noncompact manifold with bounded geometry is a manifold 
$N$ with cylindrical ends $M\times[0,\infty[$.
If we have a Clifford bundle $(E,h^E)$, also with a product structure over the cylinder, 
the Dirac operator $D$ on $E$ will in general not be Fredholm but still have finite dimensional
null space and co-null space.  
For dim$(N)=2n$ the well-known formula of Atiyah, Patodi and Singer describes
its L$^2$-index:

$$L^2\mbox{\rm-ind}(D)+\frac{h_--h_+}{2}=\int_N\widehat A(N)Ch(E/S)+
\frac{1}{2}\eta
(D^M)$$ 

Here, $\eta(D^M)$ is the Eta-invariant of the restriction $D^M$ of $D$ 
to $M.$ Die LHS is known as the {\it modified} $L^2$-index 
of $D.$

The main result shown in this paper is the analogous formula
for the ($L^2$)-$\Gamma$-index of the Dirac operator $\overline D$ on
$\overline E\to\overline N.$ 
The proof, given in Chapter 6, combines the  methods
for the proof of the $L^2$-index theorem developped in \cite{Me}, \cite{Mue2}
with the heat kernel methods used in the proof of the $\Gamma$-index theorem in 
\cite{Ro1}. The terminology is developped in Chapters 2 to 4. 

We start in Chapter 2 by 
restating the purely functional analytic description of operators on
$\Gamma$-Hilbert modules from \cite{Br}, \cite{Sh2}.
This is the natural context to develop the notion of the 
$\Gamma$-trace $\mbox{\rm tr}_\Gamma$. We give the definitions of some classes of $\Gamma$-operators
on $\Gamma$-Hilbert modules such as $\Gamma$-trace class- and $\Gamma$-compact operators,
and describe them in terms of their spectral properties.

Chapter 3 gives an outline of the analysis on manifolds of bounded geometry
as found in \cite{Sh1}, \cite{Bu1}. 
Following \cite{Ro1}, we prove the central estimates for heat kernels of elliptic
operators.
 A cornerstone of the proof of the $L^2$-$\Gamma$-index theorem will be 
the introduction 
of a spectral modification $\overline D^M+f(\overline D^M)$ of $D^M$, and we formulate the theory and the
estimates for such 'generalised elliptic' operators whenever possible.

The results of Chapters 2 and 3 will then be applied in Chapter 4 to elliptic
$\Gamma$-differential operators $\overline P$ on coverings of compact manifolds $M$.
The focus will be on a proof of the existence of the $\Gamma$-Eta-invariant
$\eta_\Gamma(\overline D^M+f(\overline D^M))$ for spectral modifications of the
Dirac operator over $M$.



Finally, Chapter 6 gives the details of the proof of the L$^2$-$\Gamma$-index theorem
\ref{th5.11} for Dirac operators over manifolds with cylindrical ends.
First, under the condition that
\begin{equation}\label{6}
0 \quad  \mbox{is an isolated point in the spectrum of}\; {\overline D}^M
\end{equation}
the analysis follows the lines of the classical case described
 in \cite{Me}. 
In the general case, a spectral modification  $\overline D_{\epsilon,u}$ of the  Dirac operator
$\overline D$ is introduced for which (\ref{6}) holds.
The main task is then to set up a good book-keeping procedure to compare the $\Gamma$-dimensions
of the null spaces of ${\overline D}_{\epsilon,u}$
and $\overline D$.


This work is the result of a project to better understand the ubiquituous
condition (\ref{6}) that appears in all generalisations of the $L^2$-index theorem. 
Using the book-keeping procedure developped in \ref{subsecl2index} the main
result would also follow from the similar result for manifolds with boundaries
given in \cite{Ra}.
In view of the $\Gamma$-signature theorem, the approach given here might perhaps
be considered to be the more natural one.

Thanks are due to Werner M\"uller who started me on this project and whose work on
Eta-invariants is at the center of many of the developments given here, as well
as to Paolo Piazza and Thomas Schick who by their kind interest revived this work.
Special thanks also to Mrs. L\"utz who reteXed the original manuscript.

\section{Operators on Hilbert $\Gamma$-modules}

This Chapter describes some of the (spectral) theory
of operators on Hilbert spaces with an equivariant action of a
discrete group  $\Gamma.$ 
Much of this material can be found in \cite{Sh2} and \cite{Br}.

\subsection{Hilbert $\Gamma$-modules and the $\Gamma$-dimension function} \label{subsect1.1}

We start with some terminology

\begin{Def}\label{d1.3}
\begin{enumerate}
\item[(a)] A free Hilbert $\Gamma$-module is a unitary right $\Gamma$-module
of the form ${\mathcal V}\otimes L^2(\Gamma),$ where ${\mathcal V}$ is a
Hilbert space.

\item[(b)] A (projective) Hilbert $\Gamma$-module is a  Hilbert space
${\H}$ with a  unitary right action of $\Gamma$, along with a 
$\Gamma$-equivariant imbedding ${\H}
\hookrightarrow{\mathcal V}\otimes L^2(\Gamma)$ into a free
Hilbert $\Gamma$-module.

\item[(c)] A morphism of Hilbert $\Gamma$-modules ${\H}_1,
{\H}_2$ is a bounded  $\Gamma$-equivariant operator
$A\in{\mathcal B}({\H}_1,{\H}_2).$ We denote the space of all such morphisms
by ${\mathcal B}_\Gamma({\H}_1,{\H}_2).$
\end{enumerate}
\end{Def}

Denote by $tr_{\mathcal V}:{\mathcal B}({\mathcal V})_+\to[0,\infty]$ the usual trace on
the Hilbert space ${\mathcal V}.$ 
We will now analyse the properties of the trace on the von Neumann Algebra
${\mathcal B}_\Gamma({\mathcal V}
\otimes L^2(\Gamma))\cong{\mathcal B}({\mathcal V})\otimes{\mathcal L}(\Gamma)$
of  endomorphisms of the free Hilbert $\Gamma$-module
 ${\mathcal V}\otimes L^2(\Gamma)$  that is induced by $tr_{\mathcal V}$ and
$tr_{\Gamma}$.

On ${\mathcal B}_\Gamma ({\mathcal V} \otimes L^2(\Gamma))$, we  have the unique f.n.s trace
 $tr_{\mathcal V}\otimes \mbox{\rm tr}_\Gamma$ which we simply denote by $\mbox{\rm tr}_\Gamma$ 
Also, for a projective Hilbert $\Gamma$-module $\mathcal H$, the $\Gamma$-embedding 
${\mathcal H}\hookrightarrow {\mathcal V}\otimes L^2(\Gamma)$
gives a $\Gamma$-embedding of operators ${\mathcal B}_\Gamma({\mathcal H} )
\hookrightarrow{\mathcal B}_\Gamma({\mathcal V} \otimes L^2(\Gamma)).$ 
Through this embedding, we can define a trace $\mbox{\rm tr}_\Gamma$  on ${\mathcal B}_\Gamma(\mathcal H)$.
It is shown in \cite{Sh2}, that this definition is in fact independent of the projective embedding.

\begin{lem}\label{l1.6}
For each (projective)Hilbert $\Gamma$-module there is a canonical f.n.s. trace 
$\mbox{\rm tr}_\Gamma$ on ${\mathcal B}_\Gamma({\mathcal H}).$ For a free Hilbert module  ${\mathcal H} ={\mathcal 
V} \otimes L^2(\Gamma)$ with orthonormal basis $(\psi_j\otimes\gamma)_
{j\in\N,\gamma\in\Gamma},$ and elements $A\in{\mathcal B}_\Gamma(\mathcal
H)_+$ this trace can be calculated by

$$\mbox{\rm tr}_\Gamma(A)=\sum_{j\in\N}\langle A\psi_j\otimes \epsilon,\psi_j
\otimes e\rangle$$
\end{lem}

\rightline{$\blacksquare$}

Using the $\mbox{\rm tr}_\Gamma$ on $\mathcal H$, a $\Gamma$-dimension function on 
$\Gamma$-invariant subspaces ${{\mathcal V}}\subset\mathcal H$ can be defined.
First, note thate the closure 
$\mbox{\rm cl}({\mathcal V})$ is $\Gamma$-invariant and denote  the orthogonal projection
onto $\mbox{\rm cl}({\mathcal V})$  by $E_{\mathcal V}\in
{\mathcal B}_\Gamma({\H}).$ The $\Gamma$-dimension of ${\mathcal V}$ is then defined by

$$\dim_\Gamma({\mathcal V}):=\mbox{\rm tr}_\Gamma([\mbox{\rm cl}({\mathcal V})])$$

The f.n.s.-property of $\mbox{\rm tr}_\Gamma$ implies the following properties for $\dim_\Gamma:$

\begin{lem}\label{l1.7}
\begin{enumerate}

\item[(a)] $\dim_\Gamma({\mathcal V})\in[0,\infty]$ and $\dim_\Gamma({\mathcal V})=0
\Leftrightarrow{\mathcal V}=\{0\}$

\item[(b)] Let $({\mathcal V}_i)_{i\in\N}$ be an increasing family of 
$\Gamma$-invariant subspaces of $\mathcal H.$ Then $\dim_\Gamma(\bigcup
{\mathcal V}_i)= \lim \dim_\Gamma({\mathcal V}_i).$

\item[(c)] Let $({\mathcal V}_i)_{i\in\N}$ be a decreasing of $\Gamma$-
invariant subspaces of $\mathcal H.$ Then $\dim_\Gamma (\bigcap{\mathcal V}
_i)
=\lim\dim_\Gamma({\mathcal V}_i).$

\item[(d)] $\dim_\Gamma$ is additive: $\dim_\Gamma({\mathcal V}_1\oplus{\mathcal V}
_2)=
\dim_\Gamma({\mathcal V}_1)+\dim_\Gamma({\mathcal V}_2).$
\end{enumerate}
\end{lem}

\rightline{$\blacksquare$}

\subsection{Classes of $\Gamma$-operators}

As in classical Hilbert space theory, we can introduce different sub-classes of the endomorphisms
 ${\mathcal B}_\Gamma({\mathcal H})$
of a Hilbert $\Gamma$-module ${\mathcal H}.$


Write  $R(A)$ for the projection onto the image  $A$
$$R(A)=\inf\{P\;|\;P\;\mbox{projection in}\; {\mathcal B}_\Gamma({\mathcal H})
\; \mbox{with}\; PA=A  \},$$
and denote the projection onto the null space of  $A$ by $N(A).$ 
From the polar decomposition of $A$ we obtain


\begin{lem}\label{l1.8} 
\begin{enumerate}
\item[(a)] Let $\sim$ denote the equivalence of projections in
in ${\mathcal B}_\Gamma({\mathcal H}).$ Then for every  
$A\in{\mathcal B}_\Gamma({\mathcal H})$

$$R(A)\sim R(A^*)=1-N(A).$$


\item[(b)]
Let $A\in{\mathcal B}_\Gamma({\mathcal H}_1,{\mathcal H}_2)$ be a 
quasiisomorphism of Hilbert $\Gamma$-modules, i.e.
 $N(A)=0$ and $R(A)={\mathcal H}_2.$ Then
$$\dim_\Gamma({\mathcal H}_1)=\dim_\Gamma({\mathcal H}_2).$$
\end{enumerate}
\end{lem}


We now define as usual

\begin{Def}\label{d1.9}
Let $\H,\H_1,\H_2$ be Hilbert $\Gamma$-modules.
\begin{enumerate}
\item[(a)] $\B_\Gamma^f(\H_1,\H_2):=\{A\in\B_\Gamma(\H_1,\H_2) \mid \mbox{\rm tr}_\Gamma(R(A))<\infty
\}$ are the $\Gamma$-operators of finite $\Gamma$-rank.

\item[(b)] The space ${{\mathcal B}}^\infty_\Gamma(\H_1,\H_2)\equiv\K_\Gamma(\H_1,\H_2)
$ 
of $\Gamma$-compact operators is the  norm closure of $\B^f_\Gamma
(\H_1,\H_2).$ 

\item[(c)] $\B^2_\Gamma(\H):=\{A\in\B_\Gamma(\H)\mid \mbox{\rm tr}_\Gamma(AA^*)<\infty\}$ are the $\Gamma$-Hilbert-Schmidt operators.

\item[(d)] $\B^1_\Gamma(\H)=\B^2_\Gamma(\H)\B^2_\Gamma(\H^*)\equiv \{A\in\B_\Gamma
(\H)\mid A= \sum^n_{i=1} S_iT^*_i$ with $S_i,T_i\in\B^2_\Gamma(\H)\}$ are
the $\Gamma$-trace class operators.
\end{enumerate}
\end{Def}

The spaces thus-defined share a number of the properties of
their classical counterparts

\begin{lem}\label{l1.10}
\begin{enumerate}
\item[(a)] $\B^f_\Gamma(\H), \B^1_\Gamma(\H), \B^2_\Gamma(\H), \B^\infty_\Gamma(\H)$ 
are two-sided *-ideals in $\B_\Gamma(\H).$

\item[(b)] $\B^1_\Gamma(\H)=\{A\in\B_\Gamma(\H)|\, \mbox{\rm tr}_\Gamma(|A|)<\infty\}.$

\item[(c)] $\B^f_\Gamma(\H)\subset \B^1_\Gamma(\H)\subset \B^2_\Gamma(\H)\subset
\B^\infty_\Gamma(\H).$

\item[(d)] $A\in\B^*_\Gamma(\H)\Leftrightarrow |A|\in\B^*_\Gamma(\H),\; *=f,1,2,
\infty.$
\end{enumerate}
\end{lem}

\begin{bem}\label{bem1.11}
Let $A,B\in\B_\Gamma(\H)$ be two self adjoint operators, and let $B$ be $\Gamma$-trace class. Then
$A$ has a unique decomposition $A=A^+-A^-$ into a sum of positive operators
$A^{\pm},\parallel A^{\pm}\parallel\le\parallel A\parallel,$ and

\begin{equation*}
\begin{split}
|\mbox{\rm tr}_\Gamma(BAB^*)|&=|\mbox{\rm tr}_\Gamma(BA^*B^*)-\mbox{\rm tr}_\Gamma(BA^-B^*)| \\
&\le\sup\{\mbox{\rm tr}_\Gamma(BA^+B^*), \mbox{\rm tr}_\Gamma(BA^-B^*)\}\le \mbox{\rm tr}_\Gamma(BB^*)\parallel
A\parallel.
\end{split}
\end{equation*}
For positive $B$ we note especially $|\mbox{\rm tr}_\Gamma(AB)|\le \mbox{\rm tr}_\Gamma(B)\parallel A
\parallel.$
\end{bem}





\subsection{The spectrum of $\Gamma$-operators}

In the following, let $\H$ be a Hilbert $\Gamma$-module, and $T:\H\supset dom(T)\to\H$ 
 a not necessarily bounded $\Gamma$-operator on $\H.$
Thus, the domain of $T$ is $\Gamma$-invariant, ${\bold r}(\gamma)dom(T)\subset dom(T)$,
and  $T{\bold r}(\gamma)\psi={\bold r}(\gamma)T\psi$ for all $\gamma\in\Gamma,\psi\in dom(T)$
( $T$ is said to be {\it affiliated} to $\B_\Gamma(\H)$). For selfadjoint $T$ the projection-valued 
measure is denoted by  $E_T(U)\in\B_\Gamma(\H)$. For each
Borel set $U\subset {\R}$ we denote the corresponding spectral subspace by 
 $\H_T(U):=\mbox{\rm im}(E_T(U))$. 
From the results of Section \ref{subsect1.1} we deduce that
$$\mu_{\Gamma,T}(U):=\mbox{\rm tr}_\Gamma(E_T(U))=\dim_\Gamma(\H_T(U))$$
defines a Borel measure on ${\R}$ whose support is the spectrum $spec(T)$ of $T$.
If  $f:{\R}\to[0,\infty]$ is a bounded Borel function, we have

$$\int_{{\R}}fd\mu_{\Gamma,T}=\mbox{\rm tr}_\Gamma(f(T)),$$

where we allow both sides of the equation to equal $\infty$.

The $\Gamma$-spectral measure $\mu_{\Gamma,T}$ gives a rough but useful classification
of the spectrum of $T.$

\begin{Def}\label{d1.13} $\mbox{\rm spec}_{\Gamma,e}(T):=\{\lambda\in{\R}\;|\;\forall
_{\epsilon>0}\;\mu_{\Gamma,T}([\lambda-\epsilon,\lambda+\epsilon])=\infty\}$ is
called the $\Gamma$-essential spectrum of $T.$
\end{Def}

This can be used to obtain the following simple spectral characterisation of $\Gamma$-compact
operators

\begin{satz}\label{s1.14}
Let $A,S$ be selfadjoint $\Gamma$-operators on $\H$. Let $S$ be bounded.
\begin{enumerate}
\item[(a)] $S\in\B^f_\Gamma(\H)\Rightarrow \mbox{\rm spec}_{\Gamma,e}(S)\subset\{0\}.$

\item[(b)] $S$ is $\Gamma$-compact $\Leftrightarrow \mbox{\rm spec}_{\Gamma,e}(S)\subset\{0\}.$

\item[(c)] Let $S$ be $\Gamma$-compact. Then ($A+S$ is selfadjoint and) 
$\mbox{\rm spec}_{\Gamma,e}(A+S)=\mbox{\rm spec}_{\Gamma,e}(A).$
\end{enumerate}
\end{satz}

{\it Proof.} (a) is obvious. (b): `$\Rightarrow$`: Choose $S$ 
compact and selfadjoint. First note that $S$ can be approximated by {\it selfadjoint} 
elements in $\B^f_\Gamma(\H)$. For $\lambda\in \mbox{\rm spec}_{\Gamma,e}(S)$
 and $\epsilon>0$ the space $\H_S(]\lambda-\epsilon,\lambda+\epsilon[)
$ is of infinite  $\Gamma$-dimension. Then choose a selfadjoint $F$ of finite 
$\Gamma$-rank, such that $\parallel S-F\parallel<\epsilon.$ For every
$\varphi\in\H_S(]\lambda-\epsilon,\lambda+\epsilon[)$ we then have

$$\parallel(F-\lambda)\varphi\parallel\le \parallel S-F\parallel
\parallel\varphi\parallel+\parallel(s-\lambda)\varphi\parallel\le 2\epsilon
\parallel\varphi\parallel,$$

thus $\H_S(]\lambda-\epsilon,\lambda+\epsilon[)\subset \H_F(]\lambda-2
\epsilon,\lambda+2\epsilon[),$ so the RHS must have infinite  $\Gamma$-dimension.
Using (a), this implies $\lambda=0.$ 

$`\Leftarrow`:$ Let $\mbox{\rm spec}_{\Gamma,e}(S)\subset\{0\}.$ Since $S$ is bounded,
the projections $E_S({\R}-]-\epsilon,\epsilon[)$ must be of finite $\Gamma$-rank for every
$\epsilon >0$. Thus
$SE_S({\R}-]-\epsilon,\epsilon[)$ 
gives a norm-approximation of $S$ by $\Gamma$-finite operators for $\epsilon \to 0$.

(c): We show $\mbox{\rm spec}_{\Gamma,e}(A)\subset \mbox{\rm spec}_{\Gamma,e}(A+S).$ To do this,
let $\lambda\in \mbox{\rm spec}_{\Gamma,e}(A),$ i.e. we have
$\dim_\Gamma(\H_A(]\lambda-\epsilon,\lambda+\epsilon[))=\infty$ for  all $\epsilon>0$.
 Now consider the set
\begin{equation*}
\begin{split}
G_\epsilon:&  =\{\varphi\in\H_A(]\lambda-\epsilon,\lambda+\epsilon[)\mid \; \parallel
S\varphi\parallel<\epsilon\parallel\varphi\parallel\}\\
 & =\H_A(]\lambda-\epsilon,
\lambda+\epsilon[)\cap \H_S(]-
\epsilon,\epsilon[)
\end{split}
\end{equation*}

But $\H_S(]-\epsilon,\epsilon[)$ is of finite $\Gamma$-codimension, and therefore
$G_\epsilon$ is of infinte  $\Gamma$-Dimension. By construction, $G_
\epsilon \subset \H_{A+S}(]\lambda-2\epsilon,\lambda+2\epsilon[)$ thus
 $\lambda\in \mbox{\rm spec}_{\Gamma,e}(A+S).$

\rightline{$\blacksquare$}

\subsection{$\Gamma$-Fredholm operators and the $\Gamma$-index}

Here we look at $\Gamma$-Fredholm operators and their properties.
Again, this closely follows the lines of 
the Hilbert space analogue. As usual, denote by $\H,\H_1,\H_2$ Hilbert
$\Gamma$-modules.

\begin{Def}\label{d1.15}
An operator $F\in\B_\Gamma(\H_1,\H_2)$ is called $\Gamma$-Fredholm if there are
 $G\in\B_\Gamma(\H_2,\H_1)$ and $K_1\in\K_\Gamma(\H_1),K_2\in\K_
\Gamma(\H_2)$, such that $$FG=1-K_2\quad GF=1-K_1.$$

Denote by ${\mathcal F}_\Gamma(\H_1,\H_2)$ the space of $\Gamma$-Fredholm operators
$\H_1\to\H_2.$
\end{Def}

From this  definition and the  ideal property of $\K_\Gamma(\H)$ we deduce
that the space of $\Gamma$-Fredholm operators $\mathcal F_\Gamma(\H)$ is closed unter the *-operation
and under concatenation of operators. It is also easy to see that $\mathcal F_\Gamma(\H)$ is an open
subset of $\B_\Gamma(\H).$

Contrary to their classical counterparts, $\Gamma$-Fredholm operators usually do not 
have a closed image, their essential spectrum can contain 0.  However, a version of
the spectral description of Fredholm operators also holds in the $\Gamma$-case.

\begin{satz}\label{s1.16} The following statements are equivalent for
a $\Gamma$-operator $F=\B_\Gamma(\H_1,\H_2)$:

\begin{enumerate}
\item[(a)] $F$ is $\Gamma$-Fredholm.
\item[(b)] $0\notin \mbox{\rm spec}_{\Gamma,e}(F^* F)$ and $0\notin \mbox{\rm spec}_{\Gamma,e}(F F^*).$
\item[(c)] $0 \notin \mbox{\rm spec}_{\Gamma,e}
\begin{pmatrix}
0  & F^* \\
F  &  0
\end{pmatrix},$ where 
$\begin{pmatrix}
0  & F^* \\
F  &  0
\end{pmatrix}\in\B_\Gamma(\H_1\oplus \H_2).$
\item[(d)] $N(F)$ is a  projection of finite $\Gamma$-rank in $\H_1$ and 
there is a projection $E$ in $\H_2$ of finite
$\Gamma$-rank, such that $\mbox{\rm im}(1-E)\subset\; \mbox{\rm im}\; (F).$ 
\end{enumerate}
\end{satz}

{\it Proof.} This follows essentially like in the classical case.

\begin{Def}\label{d1.17} Let $F\in\B_\Gamma(\H_1,\H_2)$ be $\Gamma$-Fredholm,
$E$ and as before. Then $\mbox{\rm tr}_\Gamma(N(F))<\infty$ and $\mbox{\rm tr}(1-R(F))\le tr_
\Gamma(E)<\infty.$ We can therefore define the  $\Gamma$-index of $F$ as

$$\mbox{ind}_\Gamma(F):=\mbox{\rm tr}_\Gamma(N(F))-\mbox{\rm tr}_\Gamma(1-R(F)).$$
\end{Def}
The $\Gamma$-index shares the algebraic properties of the classical index:

\begin{satz}\label{s1.18}
For operators $S,T\in\mathcal F_\Gamma(\H),K\in \K_\Gamma(\H)$ the following
holds true

\begin{enumerate}
\item[(a)] $\mbox{\rm ind}_\Gamma(S^*)=\overline{\mbox{\rm ind}_\Gamma(S)}.$
\item[(b)] $\mbox{\rm ind}_\Gamma(ST)=\mbox{\rm ind}_\Gamma(S)+\mbox{\rm ind}_\Gamma(T).$
\item[(c)] $\mbox{\rm ind}_\Gamma(S+K)=\mbox{\rm ind}_\Gamma(S)$ especially 
$\mbox{\rm ind}_\Gamma(1+K)=0.$
\item[(d)] $\mbox{\rm ind}_\Gamma:\mathcal F_\Gamma(\H)\to C$ is locally constant.
\end{enumerate}
\end{satz}

{\it Proof} The proof of these statements is analogous to the classical proofs and
can all be found in \cite{Br}.


\rightline{$\blacksquare$}

An {\it unbounded, closed} $\Gamma$-Operator $T:\H_1
\supset dom(T)\to\H_2$ is called  $\Gamma$-Fredholm, if the bounded 
 $\Gamma$-operator $T:(dom(T),\parallel\cdot\parallel_T)\to
\H_2$ is $\Gamma$-Fredholm. Here, $\parallel\cdot\parallel_T$ is
the $T$-graph norm.
We will frequently use a $\Z_2$-graded version of Proposition
\ref{s1.16}(c). An {\it unbounded, closed, odd} $\Gamma$-Operator on a $\Z_2$-graded Hilbert
$\Gamma$-module $\H=\H^+\oplus\H^-$ is called $\Gamma$-Fredholm, when 
the (unbounded) operator $T^+:\H^+\to\H^-$ is $\Gamma-$Fredholm.
The $\Gamma$-index of $T$ is then defined as $\mbox{ind}_\Gamma(T):=\mbox{ind}_\Gamma(T^+).$

\begin{satz}\label{s1.19}
Let $\H=\H^+\oplus\H^-$ be a ${\Z}_2$-graded Hilbert $\Gamma$-module,
and $T$ an {\it unbounded, closed, odd} $\Gamma$-operator on $\H.$ Then the
following two statements are equivalent

\begin{enumerate}
\item[(a)] $T$ is $\Gamma$-Fredholm.
\item[(b)] $0\notin \mbox{\rm spec}_{\Gamma,e}(T).$
\end{enumerate}
\end{satz}





\section{Manifolds of bounded geometry}

This Chapter describes some of the methods of the analysis of differential operators 
on manifolds of bounded geometry.
The results in this Chapter will be applied in the next Chapter to differential
operators on covering manifolds. Most of the concepts presented here are well-known
and can be found in different variations
in \cite{Sh1}, \cite{Bu1}, \cite{Ro1}, \cite{Lo3}, \cite{Ei}.

\subsection{Basics}

A not necessarily compact riemannian manifold $N$ of dimension
$n$ is said to be of {\it bounded geometry} if the injectivity radius $i(N)$ of $N$
is positive, and the curvature $R^N$ and all of its covariante derivatives are
bounded. A hermitian vector bundle $E\to N$ is of {\it bounded geometry},
if, in addition, the curvature  $F^E$ and all of its covariant
derivatives are bounded. Manifolds of bounded geometry admit systems of local coordinates
that have uniform $C^\infty$-estimates \cite{Ei}:

\begin{lem}\label{l2.1}
Let $E\to N$ be of bounded geometry. Fix a constant $r_0<i(N)$ and
choose a 'good' trivialisation of $E,$ i.e. $E$ is trivialised via radial parallel transport
over each ball $B(x,r_0)$.

\begin{enumerate}
\item[(a)] The metric tensor $g$ has bounded  $C^\infty$-norm w.r.t. all normal
coordinate neighborhoods $B(x,r_0)$, independent of $x$. All derivatives of coordinate change maps $\Phi_{xy}$ 
between such normal neighborhoods  $B(x,r_0)$, $B(y,r_0)$ are uniformly bounded independent of
 $x$ and $y.$
\item[(b)] In each normal neighborhood $B(x,r_0)$ with a 'good' trivialisation
of $E$, the local connection form has bounded $C^\infty$-norm independent of $x$.
All derivatives of transition maps $\Psi_{xy}$ between the trivialisations over  $B(x,r_0),$ $B(y,r_0)$ 
are uniformly bounded independent of $x$ and $y$.
\end{enumerate}
\end{lem}

\rightline{$\blacksquare$}

In the following, let  $N,E$ be of bounded geometry and choose a fixed $r_0<i(N)$.
We will always use a 'good' trivialisation of $E$ over normal neighborhoods $B(x,r_0)$.
Now define

\begin{equation*}
\begin{split}
UC^\infty(N):=\{f\in C^\infty(N)| &  \parallel\nabla^{N,k}f\parallel_\infty \;<C(k)\;
\mbox{for all}\; k\in\N\} \\
U\Gamma(N,E):=\{{\xi}\in\Gamma(N,E)|&\parallel\nabla^{E,k}\xi\parallel_\infty \;<C(k)\;\mbox
{for all}\; k\in\N\}\; \mbox{etc.}
\end{split}
\end{equation*}

Equivalently, a section $\xi$ is in $U\Gamma$ if $\xi$ and its derivatives
are uniformly bounded in any local normal coordinate neighborhood $U$ (and a corresponding
'good' trivialisation of $E$) independent of $U$.

As usual,  Sobolev spaces $H^k(N)$, $k\geq 0$, can be defined as the completion of
$C^\infty_c(N)$ with respect to the norm

$$\parallel f\parallel^2_{H^k(N)}:=\sum^k_{j=0}\parallel\nabla^jf\parallel^2_
{L^2(N,\;\otimes^jT^*N)}$$

The negative Sobolev space $H^{-k}(N)$ is then just the dual $H^k(N)$.
Similar definitions apply for spaces of sections over $N$.



Using a uniformly bounded partition of unity  (within a 'good' trivialisation of $E$ as before),
the Sobolev norms for Sections $f\in C^\infty_c(N,E)$ can locally be described as follows 

$$\parallel f\parallel^2_{H^s(N,E)}\sim\sum_{j\in N}\parallel\phi_if
\parallel^2_{H^s(U_j,\C^N)}.$$

We will use a variety of notations for the same space
$H^s(N,E)$ $\triangleq$ $H^s(N)$ $\triangleq$ $H^s(E)$ $\triangleq$ $H^s$ etc.,
depending on which part of the information is important in the particular context.

The following version of the Sobolev inequalities is now easy to prove along the
lines of its classical counterpart.

\begin{satz}\label{Sobolev} {\bf (Sobolev)}
Let $k\in\N$ and $s>k+ n/2.$ Then there is a continuous embedding
$H^s(N,E)\to UC^k(N,E).$
\end{satz}


The algebra $\mbox{\rm UDiff}^*(N,E)$ of uniform differential oprators
is generated by the uniform Sections $\Phi\in U\Gamma(N,End(E))$ and the covariant derivatives
 $\nabla^E_X$ along uniform vectorfields  $X\in U\Gamma(N,TN)$.
A differential operator $P$ is uniform if and only if its local symbol and all its derivatives
are uniformly bounded with respect to a 'good' system of coordinates of $N$ and $E$.
The operator $P\in \mbox{\rm UDiff}^k(N,E),$ maps $UC^l(N,E)$ continuously to
$UC^{l-k}(N,E)$ and maps $H^s(N,E)$ continuously to $H^{s-k}(N,E)$. Note that
the uniformity of the estimates for $P$ is essential for this to hold!
Uniform pseudodifferential operators can be defined in a similar manner.

A uniform differential operator $P$ of order  $k$ on $E$ is 
{\it uniformly elliptic,} if its principal symbol  $\sigma(P)\in U\Gamma
(T^*N,\pi^* End(E))$ has a uniform inverse outside of an $\epsilon$-neighborhood of the
null section in $T^*N$.
The construction of a parametrix for such operators can then also be performed
in a uniform manner and one can use this to show  

\begin{satz}\label{Garding} {\bf(Garding)} 
Let $T\in \mbox{\rm UDiff}^k(N,E)$ be a uniformly elliptic differential operator. Then

\begin{equation}\label{8}
\parallel\varphi\parallel_{H^{s+k}(N,E)}\le C(s,k)(\parallel\varphi\parallel_
{H^s(N,E)}+\parallel T\varphi\parallel_{H^s(N,E)})
\end{equation}

for $\varphi\in C^\infty_c(N,E),\; s\in{\R}.$
\end{satz}

\rightline{$\blacksquare$}

In Chapter 6 we will be working with spectral modifications of the Dirac operator
that are not pseudo-differential operators but which still share a number of their
mapping properties. We introduce the according spaces of operators here.
A bounded operator
$T:C^\infty_c(N,E)\to C^\infty_c(N,E)^\prime$ with Schwartz-kernel [T], 
is called an operator of order $k\in\Z,$ if it has extensions into all spaces
$\B(H^s(N,E),H^{s-k}(N,E)),s\in{\R}$. As a (possibly unbounded)
operator on $L^2(N,E)$,  $T$ is closable. To simplify things a bit we will always
ask that $T$ map components of $N_i,$ of $N$ to themselves: $T(C^\infty_c(N_i,E))\subset C^\infty_c(N_i,E)^\prime.$
Obviously, all uniform (pseudo-) differential operators of order 
$k$ are in $Op^k(N,E).$

We put the obvious family of norms onto the space $Op^k(N,E)$
and also write $Op^{-\infty}= \bigcap_{k\in N}Op^{-k}.$ 
An operator $T\in Op^k(N,E),k\ge 1,$ will be called 
{\it elliptic,} simply when it satisfies the Garding-inequality 
(\ref{8}). Note that if $T\in Op^k(N,E)$ is elliptic and $U\in Op^0(N,E),$ 
then $T+U$ is also elliptic. Also, if $T\in Op^k(N,E)$ is elliptic and 
selfadjoint, then all spectral projections of $T$ are in $Op^0(N,E).$
We note
\begin{satz}\label{s2.4} Let $T\in Op^k(N,E)$ elliptic and formally
selfadjoint, $k\ge1.$ Then $T$ is essentially selfadjoint and (without a different
notation for the closure of $T$) $dom(T)=H^k(N,E).$
\end{satz}  

\rightline{$\blacksquare$}

\subsection{Smoothing operators}\label{subsecglaettung}

Let again  $T\in Op^k(N,E)$ be {\it elliptic and formally selfadjoint.} Following
Proposition \ref{s2.4} we can interpret $T$ as a selfadjoint operator with 
$dom(T)=H^k(N,E).$ In this Section, we analyse the properties of operators
of the form $f(T)$ for a sensible choice of function $f$. The most sensible
spaces of such functions are

\begin{equation*}
\begin{split}
RB({\R}) \; := \; & \{f:{\R}\to{\C}\;\mbox{Borelfuntkion} \mid |
(1+x^2)^{k/2}f(x)|_\infty<\infty, \; k\in\N\} \\
RC({\R}) \; := \; & \{f:{\R}\to{\C}\;\mbox{stetig}\mid |
(1+x^2)^{k/2}f(x)|_\infty<\infty, \;k\in \N\}
\end{split}
\end{equation*}

The function space $RC({\R})$ with the family of seminorms $|(1+x^2)^{k/2}f(x)|_\infty$ 
is Fr\'echet. For $f\in RB({\R})$ and 
$l\in{\N}$,
the operator $T^l f(T)$ is bounded on $L^2(N,E)$ and we find using the Garding-inequality
(\ref{8})

$$\parallel f(T)\psi\parallel_{H^l(N)}  \le C(l)\sum^l_{i=0}\parallel T^if(T)
\psi\parallel_{L^2(N)}\\
\le C(l)\parallel\psi\parallel_{L^2(N)}
   \sum^l_{i=0}
|x^if|_\infty,$$
for any $\Psi \in C^\infty_c(N,E)$.
Here we have made the simplifying (but by no means essential) assumption
that $T$ is of order $1$.
Using the duality $(H^s)^*=H^{-s}$ one can show for all $k,l\in \Z,l\ge k:$

$$\parallel f(T)\psi\parallel_{H^l(N)}\le C(l,k)  \sum^{l-k}_{i=0}
\parallel T^if(T)\psi\parallel_{H^k(N)}\\
\le C(l,k)\parallel\psi\parallel_{H^k(N)}
   \sum^{l-k}_{i=0}|x^if|_\infty,$$

from which it follows that the map

$$RC({\R})\to\B(H^k(N,E),H^l(N,E)) \quad f\longmapsto f(T)$$

is continuous for all  $k,l\in\Z$, especially  $RC({\R})\to Op^{-\infty}(N,E)$
continuously. The kernels of such operators are smooth.

\begin{satz}\label{s2.5} Set $L:=[n/2+1],l\in\N.$ Then the Schwartz-kernel map

$$Op^{-2L-l}(N,E)\to UC^l(N\times N,E\boxtimes E^*)(T \longmapsto[T])$$

is continuous.
\end{satz}

{\it Proof.} We simplify notation a bit by forgetting about the coefficient bundle, i.e.
$E=\{0 \}$. Choose $r,s\in\N$ with $r+s\le l$, and an elliptic selfadjoint operator
$Q\in \mbox{\rm UDiff}^1(N)$. Then
\begin{equation*}
\begin{split}
|\nabla^r_x\nabla^s_y[T](x_0,y_0)|&\overset{\ref{Sobolev}}{\le}  C(l)\parallel\nabla^r_
x[T](x_0,\bullet)\parallel_{H^{L+s}(N)}\\
&\overset{\ref{Garding}}{\le} C(l)
\sum^{L+s}_{j=0} \parallel \nabla^r_x[TQ^j](x_0,\bullet)\parallel_{L^2(N)}
\end{split}
\end{equation*}

Setting $\xi_j(y)=\nabla^r_x[TQ^j](x_0,y)$ we can do the same estimate again

\begin{equation*}
\begin{split}
\parallel\xi_j\parallel^2_{L^2(N)}&=\int \nabla^r_x[TQ^j](x_0,y)\xi_j(y)dy =|\nabla^r_x(TQ^j\xi_j)(x_0)|\\
& \le C(l)
  \parallel TQ^j\xi_j\parallel_{H^{L+r}(N)}\\
&\le
C(l)\parallel TQ^j\parallel_{\B(L^2,H^{L+r})}
\parallel\xi_j\parallel_{L^2(N)}\\
&\le
C(l)\parallel T\parallel_
{\B(H^{-L-s},H^{L+r})}
\parallel\xi_j\parallel_{L^2(N)}
\end{split}
\end{equation*}

Together
$$|\nabla^r_x\nabla^s_y[T](x_0,y_0)|\le C(l)\sum^{L+s}_{j=0}\parallel TQ^j
\parallel_{\B(L^2,H^{L+r})}\le C(l)\parallel T\parallel_{\B(H^{-L-s},H^{L+r})}.
$$
This implies that $Op^{-2L-l}(N,E)\to UC^l(N\times E\boxtimes E^*)$ is continuous.

\rightline{$\blacksquare$}

\begin{cor}\label{c2.6} 

\begin{enumerate}
\item[(a)]  $[f(T)]\in UC^
\infty(N\times N,E\boxtimes E^*)$ for  $f\in RB({\R}).$
\item[(b)] The corresponding map $RC({\R})\to UC^\infty(N\times N,E
\boxtimes
E^*)$ is continuous.
\end{enumerate}
\end{cor}

\rightline{$\blacksquare$}

\subsection{Finite propagation speed estimates for the heat kernel}\label{subsecprop}

Sofar, for  $T$ as in Section \ref{subsecglaettung}, we have no estimates on $[f(T)]$ at infinity. 
For {\it differential operators} $P\in
\mbox{\rm UDiff}(N,E)$ which are uniformly elliptic and selfadjoint,
such estimates can be obtained by using finite propagation speed methods.
For this, we again assume that $P$ is of order $1$. From the formal selfadjointness 
of $P$ we find for any $\xi\in\Gamma_c(N,E):$

\begin{equation}\label{9}
\langle P\xi,\xi\rangle_E-\langle\xi,P\xi\rangle_E=d^*\langle\xi,
\sigma(P)\xi\rangle_E
\end{equation}

The expression $|\sigma(P)|_{2,\;T^*N\otimes End(E)}(x)$ is known as
the propagation speed of  $P$ in $x.$ The maximal propagation speed of $P$ 
on $N$ is then $c=c(P):=\sup\{|\sigma(P)|_2
(x) \mid x\in N\}:$ 
Using the spectral theorem, we know that for each $\xi_0\in C^
\infty_c(N,E)$ there is a {\it unique} solution $\xi(t)=e^{itP}\xi_
0\in L^2(N,E)$ of the wave equation $\frac{\partial\xi}{\partial t}-
iP\xi=0,$ $\xi(0)=\xi_0.$ Using \ref{Sobolev} and \ref{Garding} it is easy
to see that $\xi(t)\in UC^\infty(N,E).$ The following Lemma states, that 
$\xi(t)$ 'propagates' with finite speed:

\begin{lem}\label{l2.7} {\bf (`Energy estimate`)} For a sufficiently small $\Lambda\in
{\R}$ an all $x\in N$ the norm 
$\parallel\xi(t)\parallel_{L^2(B(x,\Lambda-ct))}$ is monotonously decreasing in $t.$
More specifically,  $P$ has propagation speed limited by $c$, since  
$\supp(\xi_0)\subset B(x,r)$ implies
$\supp(\xi(t))\subset B(x,r+ct).$
\end{lem}

{\it Proof.} This is proved using (\ref{9}) as in \cite[ 
Proposition 5.5]{Ro1}. See also the proof of Lemma \ref{l5.1}

\rightline{$\blacksquare$}

The finite propagation speed of $P$ can be used to obtain estimates for more general $f(P)$ and $[f(P)].$ 
Using the above Lemma for $\xi\in L^2(N,E)$ with $\supp(\xi)\subset B(x,r)$ we know

$$\supp(e^{isP}\xi)\subset B(x,r+c|s|).$$

This can then be plugged into the spectral representation
$$f(P)\xi=(2\pi)^{1/2}\int_{{\R}}\widehat f(s)e^{isP}\xi ds,$$
valid for all $f\in\mathcal S ({\R})$. Thus

\begin{equation}\label{f10}
\begin{split}
  \parallel f (P)\xi\parallel_{L^2(N-B(x,R))}&= \parallel(2\pi)^
{-1/2}\int_{{\R}}
\widehat f(s)e^{isP}\xi ds\parallel_{L^2(N-B(x,R))} \\
& \le\parallel(2\pi)^{-1/2}\int_{{\R}-I_R}\widehat f(s)e^{isP}\xi 
ds\parallel_{L^2(N)}\\
&\le(2\pi)^{-1/2}\parallel\xi\parallel_{L^2(N)}
 \int_{{\R}-I_R}
|\widehat f(s)|ds.
\end{split}
\end{equation}

Here, we have set $I_R:=]-\frac{R-r}{c},\frac{R-r}{c}[,$ or $I_R=\emptyset$
if $R\le r$. Thus, for  $s$ in $I_R$ the solution $e^{isP}\xi$
has not yet left $B(x,R)$. 

\begin{satz}\label{s2.8} For a sufficiently small  $r_1,$ and all 
$x,y\in N$ set $R(x,y):=\max\{0,d(x,y)-r_1\}.$ Then, writing 
$L:=[n/2+1],$ $I(x,y):=]-\frac{R(x,y)}{c},\frac{R(x,y)}{c}[$
for all $f\in \mathcal S ({\R})$, we have the estimate

$$|\nabla^l_x\nabla^k_y[f(P)](x,y)|\le C(P,l,k,r_1)\sum^{2L+l+k}_{j=0}
\int_{{{\R}}-I(x,y)}|\widehat f^{(j)}(s)|ds.$$
\end{satz}

{\it Proof.} 
Using the same technique as in 
the proof of Proposition \ref{s2.5}, we obtain





\begin{equation*}
\begin{split}
|\nabla^l_x\nabla^k_y[f(P)](x_0,y_0)| &{\le}\; C\sum^{L+l}_{i=0}
\sum^{L+k}_{j=0}
\parallel P^{j+i}f(P)\parallel_{L^2(B(x_0,r_1/2)),L^2(B(y_0,r_1/2))}\\
&\overset{\tiny(\ref{f10})}
{\le}C
\sum^{2L+l+k}_{j=0}
 \int_{{{\R}}-I(x_0,y_0)}|\widehat f^{(j)}(s)|ds.
\end{split}
\end{equation*}

\rightline{$\blacksquare$}

We now want to use this result to obtain specific estimates for the heat kernel $[f(P)]=[e^{-tP^2}]$.
It is well-known that

\begin{equation*}
\begin{split}
 \widehat f^{(k)}(s) &=\frac{1}{(2t)^{1/2}(4t)^{k/2}}
\left((4t)^{1/2}\frac{\partial}{\partial s}\right)^k
e^{-(s/(4t)^{1/2})^2} \\
&=\frac{C(k)}{t^{(k+1)/2}}H_k(s/(4t)^{1/2})
e^{-(s/(4t)^{1/2})^2},
\end{split}
\end{equation*}

where $H_k$ is the $k$th Hermite polynomial. This is even for even $k$ and odd for
odd $k$, and using
$$\int^\infty_u e^{-x^2}dx\le e^{-u^2},\quad y^s e^{-ay^2}\le \left(
\frac{s}{2ae}\right)^{s/2},\quad
s,u,y,a \in{\R}_+,$$

one obtains

$$\int^\infty_u y^s e^{-y^2}dy=\int^\infty_u y^s e^{-ey^2} e^{-(1-\epsilon)
y^2}dy\le  C(s,\epsilon) e^{-(1-\epsilon)u^2}.
$$

Using Proposition \ref{s2.8} and setting $R=R(x,y),\Lambda=2L+l+k$, we get

\begin{equation}\label{11}
\begin{split}
&  \left|\nabla^l_x\nabla^k_y [P^m e^{-tP^2}] (x,y)\right|\\
&\qquad\le C\sum^{\Lambda+n}_{j=m} t^{-j/2}\int^\infty_{R/c} 
   |H_j(s/(4t)^{1/2})|  e^{-(s/(4t)^{1/2})^2 }(4t)^{-1/2} ds \\
&\qquad \le C \sum^{\Lambda+m}_{j=m} t^{-j/2}\int^\infty_{R/2c\sqrt{t}}
|H_j(x)| e^{-x^2}dx \\
&\le C e^{-R^2/5c^2t}\sum^{\Lambda+m}_{j=m}t^{-j/2}\le
\begin{cases}
C(k,l,m,P)t^{-m/2} e^{-R^2/6c^2t},\quad t>T  \\
C(k,l,m,P) e^{-R^2/6c^2t},\; d(x,y)>2r_1,t\in{\R}_+.
\end{cases}
\end{split}
\end{equation}

This is only useful away from the diagonal.
We will obtain estimates for $[e^{-tP^2}]$ in a neighborhood of the diagonal
$N\times N$ and for small $t$ from the corresponding estimates for operators on compact manifolds.
To do this we need a 'relative' version of Proposition \ref{s2.8}.
Let $N_1, N_2$ be manifolds and $E_1\to N_1,E_2\to N_2$ 
hermitian vector bundles, all of them of bounded geometry.
On these we consider as before formally selfadjoint,  uniformly elliptic 
differential operators $P_1, P_2$  of order $1$. 

We assume that all these structures are isomorphic over an open set 
 $N_1\supset U\subset N_2$, i.e. there is a commutative diagram of
 isometries 
\begin{equation*}
\begin{CD}
N_1  \supset  
U_1  @>\phi>>U_2  \subset N_2  \\
@AAA               @AAA          \\
E_1|_{U_1}  @>\Phi>> E_2|_{U_2}
\end{CD},
\end{equation*}

such that  $P_2=\Phi P_1\Phi^{-1}$ over $U_2.$ 
Usually, we will not make these isometries explicit, but simply write
 $U\equiv U_1\equiv U_2,$ $E_1|_U\equiv E_2|_U$ etc..
Now write $c=c(P_1,P_2)$ for the maximum propagation speed of both operators
$P_1,P_2.$ We then have the following relative variant of Proposition \ref{s2.8},
see \cite{Bu1}:

\begin{satz}\label{s2.9}
Let $r_2>0.$ For  $x,y\in U$ we write $Q(x,y)=\max\{\min\{d(x,\partial U),$ $d(y,\partial U)\}
-r_2,0\}$ and $J(x,y):=]-\frac{Q(x,y)}{c},
\frac{Q(x,y)}{c}[.$ Then for $f\in \mathcal S ({\R}):$ 

\begin{equation*}
\begin{split}
|\nabla^l_x\nabla^k_y([f(P_1)](x,y)&  -[f(P_2)](x,y))| \\ 
&\le C(P_1,k,l,r_2) \sum^{2L+l+k}_{j=0}\int_{{{\R}}-J(x,y)}
|\widehat f^{(j)}(s)|ds.
\end{split}
\end{equation*}

\end{satz}

{\it Proof.} The proof is analogous to the one of Proposition \ref{s2.8}. 
The only thing to note is that the differences of operators can be estimated by





\begin{equation*}
\begin{split}
\parallel(g(P_1)&-g(P_2))\xi_j\parallel_{L^2(B(x_0,r_2/2))} \\
&=\parallel(2\pi)^{-1/2}\int_{{\R}}\widehat g(s)(e^{isP_1} -e^{isP_2})\xi_j
ds\parallel_{L^2(B(x_0,r_2 /2))} \\
& = \parallel (2\pi)^{-1/2}\int_{{\R}-J(x_0,y_0)[}\widehat g(s)(e^{isP_1}-
e^{isP_2})\xi_j ds\parallel_{L^2(B(x_0,r_2 /2))} 
\end{split}
\end{equation*}

The last equality holds, because $d(\supp(\xi_j),\partial U)\ge Q(x,y),$ 
and therefore uniqueness of the solution of the wave equation implies
$e^{isP_1}\xi_j=e^{isP_2}\xi_j$ as long as $|s|<Q(x,y)/c.$

\rightline{$\blacksquare$}

This result can now again be applied to the heat kernel, yielding in the same
manner as in (\ref{11}) for
$x,y\in U$ and $d(x,\partial U),d(y,\partial U)>r_2:$

\begin{equation}\label{12}
\begin{split}
&|\nabla^l_x\nabla^k_y\left([P^m_1e^{-tP^2_1}](x,y) -[P^m_2 e^{-tP^2_2}]
(x,y)\right)|  \\
 & \qquad\qquad\le \begin{cases}
C(k,l,m,P_1)t^{-m/2}e^{-Q(x,y)^2/6c^2t},t>T  \\
C(k,l,m,P_1)e^{-Q(x,y)^2/6c^2t},t\in{{\R}}_+
\end{cases}
\end{split}
\end{equation}

We finish this Section by adding a result on the $t \to \infty$ asymptotics
of the heat kernel. 
We know from Corollary \ref{c2.6} that the projection  $N(T):=E_T(0)$ onto the null space
of an elliptical and selfadjoint operator $T\in Op^k(N,E)$ 
has kernel $[T] \in UC^\infty$. 

\begin{satz}\label{s2.10} For $t\to\infty$ the kernel $[e^{-tT^2}]$ converges in 
$C^\infty(N\times N,E\boxtimes E^*)$ to $[N(T)].$ 
\end{satz}

{\it Proof.} This is shown in \cite[Proposition 13.14]{Ro1}.

\rightline{$\blacksquare$}


\subsection{Families of operators}

A uniform family of differential operators on an open  subset $U\subset {\R}$ 
is a uniform differential operator ${\bold P}
\overset{\wedge}{=}(P_u)_{u\in U}\in U \,\mbox{\rm Diff}^k(N\times U,E),$
which is uniformly tangential to $N$.  Locally this means that the symbols of
the uniform differential operators $P_u$ on $N$ have derivatives ( in $N$ and $U$)
of any order that can be estimated uniformly in 'good' normal neighborhoods,
independent of $u$ and the choice of neighborhood.
such uniform families can therefore be seen as $UC^\infty$-families 
$u\longmapsto P_u\in Op^k(N,E)$ with
$P^\prime_u\overset{\wedge}{=}{\bold P}^\prime(u):=[\partial_u,{\bold P}]|_u.$

For a given such family of elliptic
and selfadjoint  operators ${\bold T}:u\longmapsto T_u\in Op^k(N,E)$ 
we consider the family $e^{-t{\bold T}^2}\overset{\wedge}{=}(u\longmapsto e^{-tT^2_u})$.
Its derivative along $u$ has a well-known representation.

\begin{satz}\label{s2.12} ({\bf Duhamel formula})
 $e^{-t{\bold T}^2}$ is a differentiable map $U\subset {{\R}}$ nach $Op^{-\infty}
(N,E)$, and

$$(e^{-t{\bold T}^2})^\prime(u)=-\int^t_0e^{-sT^2_u} (T^\prime_uT_u+T_uT^
\prime_u)
e^{-(t-s)T^2_u}ds.$$
\end{satz}

Corollary \ref{c2.6} implies the differentiability of $s\longmapsto
e^{-sT^2_u}$ as a map from ${\R}_+$ to $Op^{-\infty}.$ 
To make sense of the integral at $ t \to 0$, we also need 

\begin{lem}\label{l2.13} Let $K\in Op^{-\infty}(N,E).$ Then the
maps $[0,\infty[\to Op^{-\infty}(N,E),$ $s\longmapsto Ke^{-sT^2_u},$ 
and  $s\longmapsto e^{-sT^2_u}K$ are differentiable.
\end{lem}

{\it Proof.} This is a simple exercise. Just note that
 $(1+T^2)^{-1}(e^{-sT^2}-1)$ converges in norm to $0$ for $s\to 0$.

\rightline{$\blacksquare$}

{\it Proof.} The proof of \ref{s2.12} follows the lines of the classical case by proving first continuity and
then differntiability of the family directly from the definitions.

\rightline{$\blacksquare$}

\section{Coverings of compact manifolds}

Following Atiyah's article \cite{A} we will now apply our methods 
to Hilbert $\Gamma$-modules stemming from coverings of a compact riemannian
manifold.
Section 3.3 introduces the $\Gamma$-Eta-invariant from \cite{ChG}, \cite{Ra}
and analyses its properties using the methods of \cite{Mue1}.

Start with a manifold $N$ and a vector bundle $E$ of bounded geometry.
We agree to understand by a covering of $N$ a $\Gamma$-principal bundle $\pi:\overline N\to N.$ 
This is what we mean when referring to $\Gamma$ as the 'covering group'.

Fix a fundamental domain  ${\mathcal F}\subset\overline N.$ For
$L \subset N,$ we write the lift as $\overline L:=\pi^{-1}(L)$ qnd define ${\mathcal F}(L)
$ to be
${\mathcal F}\bigcap \overline L.$ The lift of the vector bundle  $E$ to
$\overline N$ is denoted by  $\overline E,$ and  $\overline\xi:
\overline N\to\overline E$ is meant to refer to the lift of the section
$\xi:N\to E$ etc..
From the bounded geometry of  $N,E$ it follows immediately that 
 $\overline N,\overline E$ also have bounded geometry.

\subsection{Coverings and Hilbert $\Gamma$-modules}

We show that the Sobolev spaces over  $\overline N$ can be interpreted as Hilbert
$\Gamma$-modules. First, note that

$$L^2(\overline N,\overline E) \overset{\sim}{\to} L^2({\mathcal F},
\overline E|
_{{\mathcal F}})\otimes L^2(\Gamma)
\left( \xi\longmapsto\sum_{\gamma\in\Gamma}\xi \gamma|_{{\mathcal F}}
\otimes
\gamma^{-1}\right)$$

is an isomorphism of (right) $\Gamma$-modules (i.e.  $(\xi\gamma)(x)=\xi(x\gamma^{-1}).$ )
Thus, $L^2(\overline N,\overline E)$ is a 
free Hilbert $\Gamma$-module.

For Sobolev spaces other than $L^2$ the above map is not usually
an isomorphism. However, for $s>0$ the map

$$H^s(\overline N,\overline E)\to H^s({{\mathcal F}},\overline E |_{
{\mathcal F}})\otimes L^2(\Gamma)$$

is still an isometric imbedding, which is the requirement for $H^s(\overline N,\overline E)$ 
to be a Hilbert $\Gamma$-module.

We now have the following generalisation of Rellich's Lemma:

\begin{satz}\label{s3.1} ({\bf Rellich}) Let $f\in C^ \infty_c(N)$ a compactly
supported function. For $s,s^\prime\in\R$ and $s>s^\prime$ the map

$$H^s(\overline N,\overline E)\overset{M_{\overline f}}{\to} H^s
(\overline N,\overline E) \overset{\iota}{\hookrightarrow}H^{s^\prime}
(\overline N,\overline E)$$

is a  $\Gamma$-compact morphism of Hilbert $\Gamma$-modules.
\end{satz}

{\it Proof.} Instead of $M_{\overline f}$ we write $\overline f,$ and we drop the
explicit mention of the embedding $\iota$. Obviously, the maps
$\overline f,\iota$ are $\Gamma$-operators. Now for $s,s^\prime>0$, we have  
a commutative diagram:

\begin{equation*}
\begin{CD}
H^s({\overline N},{\overline E})@>{{\overline f}}>>H^{s^\prime}({\overline N},
{\overline E}) \\
  @VVV     @VVV\\
H^s({\mathcal F},{\overline E}|_{\mathcal F})\otimes L^2(\Gamma)@>{{\overline f}
|_{\mathcal F}\otimes I}>>H^{s^\prime}({\mathcal F},{\overline E}|_{\mathcal F})
\otimes L^2(\Gamma)
\end{CD}.
\end{equation*}

But $\supp(\overline f)\cap {\mathcal F}$ is compact and the map
$H^s({\mathcal F},\overline E|_{\mathcal F})\overset{\overline f|_{\mathcal F}}
{\to}H^{s^\prime}({{\mathcal F}},\overline E|_{\mathcal F})$ 
is compact in the classical sense due to the classical Rellich-Lemma.
Thus,  $\overline f|_{\mathcal F}\otimes I$ is $\Gamma$-
compact and the assertion follows from the fact that the  vertical arrows
in the diagram are isometric embeddings.

\rightline{$\blacksquare$}

We will now look into the description of the Schwartz-kernels of the $\Gamma$-trace class and 
$\Gamma$-Hilbert-Schmidt operators
on $L^2(\overline N,\overline E)$.

\begin{satz}\label{s3.2} Let $A\in\B_\Gamma(L^2(\overline N,\overline E)).$ Then

\begin{enumerate}
\item[(a)] $A\in\B^1_\Gamma(L^2(\overline N,\overline E))\Leftrightarrow
\chi_{\mathcal F}|A|\chi_{\mathcal F}\in\B^1(L^2({\mathcal F},\overline E|_
{\mathcal F})).$
\item[(b)] $A\in\B^1_\Gamma(L^2(\overline N,\overline E))\Rightarrow
\mbox{\rm tr}_\Gamma(A)=\mbox{\rm tr}(\chi_{\mathcal F} A\chi_{\mathcal F}).$ If the kernel of  $A$ 
is continuous:
$$\mbox{\rm tr}_\Gamma(A)=\int_{\mathcal F}tr_{\overline E}([A])(x,x)d\mbox{\rm vol}_{\mathcal F}
(x)=\int_N\pi_*tr_{E}([A])(x,x)d\mbox{\rm vol}_N(x),$$
where we have used $[A](x\gamma,x\gamma)=[A](x,x)$ to push down
$[A](x,x)$ on the diagonal to  $\pi_*([A](x,x))$ on the basis $N$.
\item[(c)] $A\in\B^2_\Gamma(L^2(\overline N,\overline E))\Leftrightarrow
[\chi_{\mathcal F}A]\in L^2({\mathcal F}\times\overline N,\overline E
|_{\mathcal F}\boxtimes \overline E^*).$
\end{enumerate}
\end{satz}

{\it Proof.} (a) and (b): Let $(\psi_j)_{j\in\N}$ be
an orthonormal basis of $L^2({\mathcal F},{\overline E}_{\mathcal F})$ $\subset
L^2(\overline N,\overline E).$
Then $(\psi_j\gamma)_{j\in\N,\gamma\in\Gamma}\cong(\psi_j\otimes
\gamma)_{j\in\N,\gamma\in\Gamma}$ is an orhtonormal basis of
$LÖ^2(\overline N,\overline E)\cong L^2({\mathcal F},\overline E|_{\mathcal F})
\otimes L^2(\Gamma)$ and Lemma \ref{l1.6} implies:

\begin{equation*}
\begin{split}
\mbox{\rm tr}_\Gamma(|A|)\overset{\wedge}{=}(\mbox{\rm tr}_\Gamma\otimes tr)(|A|)&  =\sum_{j\in
\N}\langle |A|\psi_j\otimes e,\psi_j\otimes e\rangle\\
&=\sum_{j\in\N}
\langle|A|\psi_j,\psi_j\rangle
  =\mbox{\rm tr}(\chi_{\mathcal F}|A|\chi_{\mathcal F}).
\end{split}
\end{equation*}

Thus, the operator $A$ is $\Gamma$-trace class, if and only if $\chi_{\mathcal
F}A_{\chi{\mathcal F}}$ is trace class. The integral representation of $A$ directly follows.

(c): $AA^*$ is positive, thus

\begin{equation*}
\begin{split}
AA^*\in\B^1_\Gamma(L^2(\overline N,\overline
E))  &   \overset{(a)}{\Leftrightarrow}
\chi_{\mathcal F}AA^*\chi_{\mathcal F}\in\B^1(L^2({\mathcal
 F},\overline E|_{\mathcal F}))  \\
 &  \Leftrightarrow
[\chi_{\mathcal F} A]\in L^2({\mathcal F}\times \overline N,\overline E|_
{\mathcal F}\boxtimes \overline E^*).
\end{split}
\end{equation*}

\rightline{$\blacksquare$}


\subsection{Elliptic operators}

Let now $M$ be a {\it compact} riemannian manifold and $E$ a 
$\Z_2$-graded vector bundle over $M.$
Denote by $Op^k_\Gamma(\overline M,\overline E)$, $k>0$, the subspace of  $\Gamma$-equivariant operators
 in $Op^k(\overline M,\overline E).$ Of course, the typical example will be the lifts 
$\overline P\in \mbox{\rm UDiff}^k(\overline M,\overline E)$ of differential operators
$P \in \mbox{\rm Diff}^k(M,E).$

From Proposition \ref{Garding}, for an operator $T\in Op^k_\Gamma(\overline M,\overline E)$ that is elliptic and selfadjoint,
the operator $(T\pm i)^{-1}$ is in $\B_\Gamma(H^s(\overline M,
\overline E),H^{s+k}(\overline M,\overline E)),$ and thus $\Gamma$-compact according to Proposition
\ref{s3.1}:
$$(T\pm i)^{-1}\in\K_\Gamma(L^2(\overline M,\overline E).$$
Especially $\mbox{spec}_{\Gamma,e}(T^2+1)\subset\{0\},$ i.e. $T$ is $\Gamma$-
Fredholm.

More generally, set
$n=\dim(M),$ $L=[n/2+1]$  and write$Op^m_\Gamma(\overline M,\overline E)_+$ 
for the elements in  $Op^m_\Gamma(\overline M,\overline E),$  that are positive operators on
$L^2(\overline M,\overline E).$

\begin{satz}\label{s3.4} 
\begin{enumerate}
\item[(a)] The elements of $Op^{-2L}_\Gamma(\overline M,
\overline E)$ are $\Gamma$-trace class, and the map $Op^{-2L}_\Gamma
(\overline M,\overline E)_+\to \B^1_\Gamma(\overline M,
\overline E)$ is continuous.
\item[(b)] The elements in $Op^{-L}_\Gamma(\overline M,\overline E)$ are $\Gamma$-Hilbert-Schmidt, 
and the map $Op^{-L}_\Gamma(\overline M,\overline E)
\to \B^2_\Gamma(L^2(\overline M, \overline E))$ is continuous.
\item[(c)] $RC({\R})\to\B^1_\Gamma(L^2(\overline M,\overline E))(f\longmapsto
f(T))$ is continuous.
\end{enumerate}
\end{satz}

{\it Proof.} (a): From Proposition \ref{s2.5} the operator $A\in Op^{-2L}_\Gamma(\overline 
M,\overline E)$ has uniformly continuous Schwartz-kernel $[A]$, which therefore can be integrated
over ${\mathcal F}\subset \Delta\subset \overline M\times
\overline M.$ Proposition \ref{s3.2} then implies that $A$ is $\Gamma$-trace class, 
if $A$ is a {\it positive} operator. The continuity of the kernel-map then follwos from the estimate
$$\mbox{\rm tr}_\Gamma(|A|)=\int_{\mathcal F}[A](x,x)d\mbox{\rm vol}_{\overline M}(x)\le
|[A]|_\infty vol ({\mathcal F}).$$
This proves the second part of (a). Part (b) follows from the continuity of
$$Op^{-L}_\Gamma(\overline M,\overline E)\to Op^{-2L}_\Gamma(\overline M,
\overline E)_+\quad (A\longmapsto AA^*)$$
Now, for the first part of (a) choose an  elliptic
differential operator $T\in \mbox{\rm Diff}^L(\overline M,E)$ such that $ST=1-R,$ with 
suitable parametrix $S\in Op^{-L}_\Gamma$ and error term $R\in Op^{-\infty}_\Gamma.$ 
Given a (not necessarily positive) $A\in Op^{-2L}_\Gamma$, part (b) then implies
that the operators $S,TA,R,A\in\B^2_\Gamma,$ thus $A=(ST+R)A=S(TA)+RA\in\B^1_
\Gamma.$
(c) now follows from (a) and the continuity of
$$RC({{\R}})\overset{|\cdot|}{\to}RC({\R})_+\to Op^{-2L}_\Gamma(\overline
M,\overline E)_+.$$

\rightline{$\blacksquare$}

\begin{cor}\label{c3.5}

As before, let $T\in Op^k_\Gamma(\overline M,\overline E)$ elliptic and selfadjoint.
Then $\mbox{spec}_{\Gamma,e}(T)=\emptyset$ and the spectral measure $\mu_{\Gamma,T}$ is polynomially bounded. 
More precisely, the spectral 'counting function' satisfies the estimate $N_{\Gamma,T}(\lambda):=\dim_\Gamma(\H_T(]-
\lambda,\lambda[))\le C \lambda^{2L/k}.$
\end{cor}

{\it Proof.} The operator $(1+T^2)^{-L/k}$ is $\Gamma$-trace class, thus
$$\int_{\R}(1+x^2)^{-L/k}d\mu_{\Gamma,T}(x)=\mbox{\rm tr}_\Gamma((1+T^2)^{-L/k})
<\infty.$$

and we can write for the spectral counting function

$$N_{\Gamma,T}(\lambda)=\int_{\R}\chi_{]-\lambda,\lambda[}(x)d\mu_{\Gamma,T}(x)\le(1+\lambda
^2)^{L/k}\int_{\R}(1+x^2)^{-L/k}d\mu_{\Gamma,T}(x).$$

\rightline{$\blacksquare$}

In the case that  $E\rightarrow N$ is $\Z_2$-graded, $T$ odd, we can use these
results to calculate the index:

\begin{lem}\label{l3.6} $ \mbox{\rm ind}_\Gamma(T)=\mbox{\rm str}_\Gamma(e^{-tT^2}).$
\end{lem}

{\it Proof.} The proof now runs parralel to the proof for compact manifolds.
The map $f\mapsto f(T)$ is continuos from
$RC(\R)$ to $\B^1_\Gamma(L^2(\overline M,\overline E)).$
The family $t\mapsto e^{-tx^2}$ is differentiable as a map from
$\R_+$ to $RC(\R).$ Therefore  $\mbox{\rm str}_\Gamma(e^{-tT^2})$ 
is differentiable in $t>0$ and 

$$\frac{d}{dt}\mbox{\rm str}_\Gamma(e^{-tT^2})=-\mbox{\rm str}_\Gamma(T^2 e^{-tT^2})=-\frac{1}{2}
\mbox{\rm str}_\Gamma([T,Te^{-tT^2}])=0.$$

This shows that $\mbox{\rm str}_\Gamma(e^{-tT^2})$ is independent of $t$. 
But according to Proposition \ref{s2.10}, the heat kernel $[e^{-tT^2}]$ converges 
 in $C^\infty$ to $[N(T)]$ for $t \to \infty$. Thus 

\begin{equation*}
\begin{split}
\mbox{\rm str}_\Gamma(N(T))&  =\int_{\mathcal F}\mbox{\rm str}^{\overline E}[N(T)](x,x) d\mbox{\rm vol}_{
\overline M}(x)\\
&  =\lim_{t\to\infty}\int_{\mathcal F}\mbox{\rm str}^{\overline E}
[e^{-tT^2}](x,x)d\mbox{\rm vol}_{\overline M}(x)=\mbox{\rm str}_{\Gamma}(e^{-tT^2})
\end{split}
\end{equation*}

where the LHS is of course $\mbox{ind}_\Gamma(T).$

\rightline{$\blacksquare$}

In the special case that $T=\overline P$ is the lift of an odd operator
$P\in \mbox{\rm Diff}^1(M,E)$, the operator $ \overline P$ can locally be compared to  $P$.
For small $\epsilon>0$ and any $\overline x\in\overline M$, the
ball  $B(\overline x,\epsilon)\subset \overline M$ is isometric to the ball
 $B(x,\epsilon)\subset M$, $x=\pi(\overline x)$, and the vector bundles and operators
are isometric as well, i.e.
$\overline E|_{B(\overline x,\epsilon)} \backsimeq E|_{B(x,\epsilon)}$ and
$\overline P|_{B(\overline x,\epsilon)} \backsimeq P|_{B(x,\epsilon)}.$ 
Applying the estimate (\ref{12}) then gives for  
$t\in{\R}_+$ and suitable constants $C_1,c_2>0$

\begin{equation}\label{14}
|[e^{-t\overline P^2}](\overline x,\overline x)-[e^{-tP^2}](x,x)|\le C_1 
e^{-c_2/t}.
\end{equation}

This allows us to reprove Atiyah's theorem
\begin{theo}\label{th3.7} {\bf (Atiyah)} $\mbox{\rm ind}_\Gamma(\overline P)=\mbox{\rm ind}(P).$
\end{theo}

{\it Proof.} We just need to put together the results obtained thusfar

\begin{equation*}
\begin{split}
\mbox{ind}_\Gamma({\overline P}) & = \lim_{t\to 0}\mbox{\rm str}_\Gamma(e^{-t\overline P^2})
=\lim_{t\to 0}\int_{\mathcal F}\mbox{\rm str}^{\overline E}[e^{-t\overline P^2}]
(\overline x,\overline x)d\mbox{\rm vol}_{\overline M}(\overline x)  \\
& \overset{\ref{14}}{=}\lim_{t\to 0}\int_M \mbox{\rm str}^E[e^{-tP^2}](x,x)d\mbox{\rm vol}_M(x)=\lim_
{t\to 0}\mbox{\rm str}(e^{-tP^2})=\mbox{ind}(P).
\end{split}
\end{equation*}

\rightline{$\blacksquare$}

\subsection{The $\Gamma$-Eta-invariant}\label{subseceta}

Let $T\in Op^1_\Gamma(\overline M,\overline E)$ elliptic and 
selfadjoint.
 We will try to define the $\Gamma$-Eta-invariant for $T$ as
the value in $s=0$ of




\begin{equation}\label{16}
\eta(T)(s)=\frac{1}{\Gamma(\frac{s+1}{2})}\int^\infty_0 t^{(s-1)/2}
\mbox{\rm tr}_\Gamma(T e^{-tT^2})dt.
\end{equation}

This will require an analysis of the expression $\mbox{\rm tr}_\Gamma(Te^{-tT^2})$ for
$t\to\infty$ and $t\to0$.

The first part of the analysis can be done for  $T=\overline P,$ with an elliptic and selfadjoint 
differential operator $P\in \mbox{\rm Diff}^1(M,E).$ 
For $t\to 0$, there is the well-known asymptotic development on $M:$

$$\mbox{\rm tr}(Pe^{-tP^2})\sim\sum^\infty_{j=0}b_j(P)t^{(j-n-1)/2}.$$

The coefficients $b_j(P)$ are local integrals $b_j(P)=\int_M\beta_j(P),$ where $\beta_j(P)(x)$ 
only depends on the symbol of $P$ and its derivatives. As in the proof of Theorem \ref{th3.7},
the estimate (\ref{14}) implies that $\mbox{\rm tr}_\Gamma(\overline P
e^{-t\overline P^2})$ and $\mbox{\rm tr}(Pe^{-tP^2})$ have the same asymptotics for $t \to 0$

\begin{equation}\label{17}
\mbox{\rm tr}_\Gamma(\overline P e^{-t\overline P^2})\sim\sum^\infty_{j=0}b_j
(\overline P)t^{(j-n-1)/2},
\end{equation}

and $b_j({\overline P})=b_j(P)=\int_{\mathcal F}\beta_j(\overline P).$
Thus the expression
$$\eta_\Gamma({\overline P})(s)_\kappa:=\frac{1}{\Gamma(\frac{s+1}{2})}\int^
\kappa_0t^{(s-1)/2}\mbox{\rm tr}_\Gamma(\overline P e^{-t\overline P^2})dt$$

exists for $\kappa\in{\R}_+,s>n+1,$ and has an asymptotic development of the form

\begin{equation*}
\begin{split}
\eta_\Gamma(\overline P)(s)_\kappa &  \sim \frac{1}{\Gamma(\frac{s+1}{2})}\sum
^\infty_{j=0}b_j({\overline P})
\frac{\kappa^{(j-n)/2+s/2}}{(j-n)/2+s/2}.
\end{split}
\end{equation*}

The expression $\Gamma((s+1)/2)\eta_\Gamma(\overline P)(s)_\kappa$ therefore
has a meromorphic continuation to  $\C$ with singularities of $1$st order (at most)
in $s\in\{n-j|j\in\N\}.$ Thus, $\eta_\Gamma(\overline P)(s)_\kappa$ is
holomorphic in $s=0,$ exactly when $b_n(\overline P)=0.$ An equivalent formulation is
to say that the development

$$\int^\kappa_\delta t^{-1/2}\mbox{\rm tr}_\Gamma(\overline P e^{-t\overline P^2})dt=
\sum_{j\neq n}b_j(\overline P)\frac{2}{j-n}(\kappa^{\frac{j-n}{2}}-
\delta^{\frac{j-n}{2}})$$

in powers of $\delta$ exists. For a function $f(\delta)$ permitting such
a development in $\delta \to 0$ we denote the constant coefficient in the development by
 LIM$_{\delta\to 0}f(\delta)$.
Thus, for $b_n(\overline P)=0$ 

$$\eta_\Gamma(\overline P)(0)_\kappa=  \frac{1}{\Gamma(\frac{1}{2})}
\sum_{j\neq n}b_j(\overline P)\frac{2}{j-n}\kappa^{\frac{j-n}{2}} =
\mbox{LIM}_{\delta\to 0}\int^\kappa_\delta
\frac{t^{-1/2}}{\Gamma(\frac{1}{2})}(\overline P e^{-t\overline P^2})dt.$$

In the case of the Dirac operator  $\overline D$, the coefficients $b_j(\overline D)$ 
vanish for  all $j\le n$ and $\eta_\Gamma(\overline D)(s)_\kappa$ is
holomorphic in $0,$ i.e. it exists without any regularisation. 

The analysis of the term

$$\eta_\Gamma(\overline P)(s)^\kappa:=\frac{1}{\Gamma(\frac{s+1}{2})}
\int^\infty_\kappa t^{(s-1)/2}\mbox{\rm tr}_\Gamma(\overline P e^{-t\overline P^2})dt$$

is easy, because this integral exists for $s\le0$ due to the polynomial boundedness of
$\mu_{\Gamma,\overline P}$.
%


For the Dirac operator, we can thus define:

$$\eta_\Gamma(\overline D)=\eta_\Gamma(\overline D)(0)_\kappa+\eta_\Gamma(
\overline D)(0)^\kappa.$$

We have to extend this definition of the $\Gamma$-Eta-invariant to  modifications  $Q=\overline D +K$ 
of $\overline D$ with  $K \in Op^{-\infty}(\overline M)$.
The main example we have in mind are modifications of $\overline D$
where $K=f(\overline D )$ for a bounded measurable function $f:{\R}\to{\R}$. First

\begin{equation*}
\begin{split}
Qe^{-tQ^2}&-\overline D e^{-t\overline D^2}\\
&=K e^{-t(\overline D + K)^2}-\overline D \int_0^t e^{-s(\overline D + K)^2}( K\overline D+ \overline D K +K^2)
 e^{-(t-s){\overline D}^2}ds.
\end{split}
\end{equation*}

 Lemma $\ref{l2.13}$
implies that this converges for $t\to 0$ to $K$  in  $Op^{-\infty}(\overline M)$. 
This gives an asymptotic development

\begin{equation}\label{19}
\mbox{\rm tr}_\Gamma(Qe^{-tQ^2})\sim \sum^{n+1}_{j=0}b_j(\overline D)t^{(j-n-1)/2}+tr_
\Gamma(K)+g(t),\quad t\to 0,
\end{equation}

where $g$ is a continuous on $[0,\infty[$ and $g(0)=0$. 
Since  $b_j(\overline D)=0$ for $j\le n$, the expression $\eta_\Gamma(Q)(0)_\kappa$ 
exists.Using again the polynomial boundedness of the measure $\mu_{\Gamma,\overline Q}$.
\par
\medskip
To compare the $\Gamma-\eta$-invariants for $\overline D$ and ${\overline D} + K$, consider now the 
differentiable family $u\mapsto Q_u:=\overline
D+K+u$ of elliptic, selfadjoint operators in 
$Op^1_\Gamma(\overline M,\overline E).$ From Proposition \ref{s2.12}, the trace
$\mbox{\rm tr}_\Gamma(Q_ue^{-tQ^2_u})$ is differentiable  in $u$ and using the
trace property one has

\begin{equation*}
\begin{split}
\frac{\partial}{\partial u}\mbox{\rm tr}_\Gamma(Q_ue^{-tQ^2_u})&  =\mbox{\rm tr}_\Gamma\left(Q^\prime
_ue^{-tQ^2_u}-t(Q^\prime_uQ_u+Q_uQ^\prime_u)e^{-tQ^2_u}\right)\\
&  =
\left( 1+2t \frac{\partial}{\partial t}\right)\mbox{\rm tr}_\Gamma
(Q^\prime_ue^{-tQ^2u}).
\end{split}
\end{equation*}

From the asymptotic development of $\mbox{\rm tr}_\Gamma(Q_u e^{-tQ_u^2}),$ one deduces that $\eta_\Gamma(Q_u)(s)_\kappa$ 
exists for large $s$ and can be continued meromorphically. 
Partial integration gives

\begin{equation*}
\begin{split}
\frac{\partial}{\partial u}&\eta_\Gamma(Q_u)(s)_\kappa   =\frac{\partial}
{\partial u}\int^\kappa_0 \frac{t^{(s-1)/2}}{\Gamma(\frac{s+1}{2})}
\mbox{\rm tr}_\Gamma(Q_u  e^{-tQ^2_u})dt  \\
&  = \int^\kappa_0\frac{t^{(s-1)/2}}{\Gamma(\frac{s+1}{2})}(1+2t\frac{\partial}
{\partial t})\mbox{\rm tr}_\Gamma(Q^\prime_u e^{-tQ^2_u})dt  \\
& = \frac{2\kappa^{(s+1)/2}}{\Gamma(\frac{s+1}{2})} \mbox{\rm tr}_\Gamma(Q^\prime_u
e^{-\kappa Q^2_u})-\frac{s}{\Gamma(\frac{s+1}{2})}
\int^\kappa_0 t^{(s-1)/2}\mbox{\rm tr}_\Gamma(Q^\prime_u e^{-tQ^2_u})dt.
\end{split}
\end{equation*}
Here, we have used the identity:
$$ t^{(s-1)/2}(1+2t\frac{\partial}
{\partial t})h = (2\frac{\partial}
{\partial t}t -s) t^{(s-1)/2}h.$$

In order to understand the asymptotic development of $\frac{\partial}{\partial u}\eta_\Gamma(Q_u)(s)_\kappa$,
we need to understand the asymptotic development of the last integral.
In the case at hand, $Q^\prime_u$ is the identity, and
the asymptotic development of the heat kernel for $\overline D +u$
gives for $t\rightarrow 0$:
$$\mbox{\rm tr}_\Gamma(Q^\prime_ue^{-tQ^2_u})= \mbox{\rm tr}_\Gamma(e^{-t(\overline D + K +u)^2})\sim\sum^n_{j=0}a_j(\overline D+u)t^{(j-n)/2}
+g(t),$$
where again $g$ is continuous on $[0,\infty[$ ith $g(0)=0$. 
The integral therefore has a meromorphic extension to a neighborhood of $0$ in ${\C}$
with a singularity of order at most $1$ in $0$.
But then, $\frac{\partial}{\partial u}\eta_\Gamma(Q_u)(s)_\kappa$ 
is holomorphic  in $0$ and 

$$\frac{\partial}{\partial u} Res_{s=0} \;\eta_\Gamma(Q_u)(s)_\kappa=Res_{s=0}
\frac{\partial}{\partial u}\eta_\Gamma(Q_u)(s)_\kappa=0,$$

and  $Res_{s=0}\eta_\Gamma(Q_u)(s)_\kappa$ is constant in $u.$ Since $\eta_
\Gamma(Q_0)(s)_\kappa$ is holomorphic in $0$ so is
$\eta_\Gamma(Q_u)(s)_\kappa$ holomorphic in $0.$

In Chapter 6 we will look at families $Q_u=\overline D+u-\Pi\overline D,$ with
$\Pi= E_{\overline D}(]-\epsilon,\epsilon[)$. For this special case we note

\begin{lem}\label{l3.8}
For $Q_u=\overline D+u-\Pi\overline D$  the $\Gamma$-Eta-invariant 

$$\eta_\Gamma(Q_u)=LIM_{\delta\to0}\int^\kappa_\delta\frac{t^{-1/2}}{\Gamma
(\frac{1}{2})}\mbox{\rm tr}_\Gamma(Q_ue^{-tQ^2_u})dt+\int^\infty_\kappa
\frac{t^{-1/2}}{\Gamma(\frac{1}{2})}\mbox{\rm tr}_\Gamma(Q_ue^{-tQ^2_u})dt$$

exists for all $u,$ and for $\epsilon>0$ we have

\begin{enumerate}
\item[(a)] $\eta_\Gamma(Q_u)-\eta_\Gamma(Q_0)=\mbox{\rm sgn}(u)\mbox{\rm tr}_\Gamma(\Pi)$ also
$\eta_\Gamma(Q_0)=\frac{1}{2}(\eta_\Gamma(Q_u)+\eta_\Gamma(Q_{-u})).$
\item[(b)] $|\eta_\Gamma(\overline D)-\eta_\Gamma(Q_0)|=|\eta_\Gamma(\Pi
\overline D)|\le \mu_{\Gamma,\overline D}(]-\epsilon,\epsilon[-\{0\}).$
\end{enumerate}
\end{lem}

\rightline{$\blacksquare$}

\subsection{Residually finite coverings}

As usual, let $M$ be a compact manifold, $\pi:\overline M\to M$ a
$\Gamma$-principal bundle.  We assume that the the covering group $\Gamma$ is {\it finitely generated,}
and fix a word metric $|\bullet|.$
The following is a standard result in geometric
group theory.

\begin{lem}\label{l3.9} Let $\overline M$ be connected. Then
$\overline M$ and $\Gamma$ are quasiisometric, i.e. there are constants
$A,B>0$ with

$$B^{-1}|\gamma|-A\le d(x,x\gamma)\le B|\gamma|+A$$

for all $x\in\overline M,$ $\gamma\in\Gamma.$ 
\end{lem}

\rightline{$\blacksquare$}

The considerations in this Section are based on the following Proposition
that shows how the heat kernel of a selfadjoint elliptic differential operator
 $P\in \mbox{\rm Diff}^1(M,E)$ can be reconstructed from
the heat kernel of its lift $\overline P$.

\begin{satz}\label{s3.10} Choose $\overline x$ and write $x=\pi(\overline x)$.
Then

$$[P^l e^{-tP^2}](x,y)=\sum_{\gamma\in\Gamma}[\overline P^l e^{-t\overline P}
](\overline x,\overline y\gamma),$$

and the RHS converges in $UC^\infty$ uniformly in $0<a\le t\le b\le\infty$ .
\end{satz}

{\it Proof.} The proof can be found in \cite{Lo3}.
Here, we only show the absolute convergence of the RHS independent of $\overline x,\overline y$ 
and $t\in[a,b]$.
It suffices to show this for $\overline x,\overline y\in{\mathcal F}$:

\begin{equation*}
\begin{split}
& \sum_{\gamma\in\Gamma}|[\overline P^l e^{-t\overline P^2}](\overline x,
\overline y\gamma)| \overset{\tiny(\ref{11})}{\le} C+C\sum_{e\neq\gamma\in\Gamma}
t^{-l/2}e^{-\delta(d(\overline x,\overline y\gamma)-r_0)_+^{2}/t}  \\
& \le C+C(a,b)\sum_{e\neq \gamma\in  \Gamma}e^{\delta d(\overline x,
\overline xy)^2/b}\le C+C(a,b)\sum_{e\neq\gamma\in\Gamma}
e^{-\delta B|\gamma|^2/b}<\infty
\end{split}
\end{equation*}

and similarly for the derivatives. 

\rightline{$\blacksquare$}

Let now $\overline M\to M$ be {\it residually finite}, i.e.
there is a tower of groups

$$\Gamma=\Gamma_0\vartriangleright \Gamma_1 \vartriangleright \ldots
\vartriangleright \Gamma_i \vartriangleright \Gamma_{i+1}\ldots\Gamma_\infty=
\{e\},$$

i.e. $\Gamma_{i+1}$ is normal of finite index in $\Gamma_i.$ 
Setting $M_i:=\overline M/\Gamma_i$ the covering $M_i\to M$ is finite
with deck-transformations $\Gamma/ \Gamma_i.$ If $\overline M\to M$ is
residually finite then so is $\overline M^0\to M$.
Let now $E$ be a over $M$ with Dirac operator
$D$ and denote the lifts of the vector bundle $E$ to $M_i$ by $E_i$ and the lifted Dirac operators by $D_i.$ 
The following is a simple consequence of Proposition \ref{s3.10}

\begin{lem}\label{l3.11}
Write $b_i:=\dim(\ker(D_i)),$ $b_\Gamma:=\dim_\Gamma(\ker(\overline D)),$
$d_i:=[\Gamma:\Gamma_i].$ Then

\begin{enumerate}
\item[(a)] $\lim_{i\to\infty}\frac{1}{d_i}\mbox{\rm tr}(e^{-tD^2_i})=\mbox{\rm tr}_\Gamma(e^{-t
\overline D^2})$
\item[(b)] $\lim \sup_{i\to\infty}\frac{b_i}{d_i}\le b_\Gamma$
\end{enumerate}
\end{lem}

\rightline{$\blacksquare$}

Choose an increasing sequence $\{e\}=S_0\subset S_1\subset S_2\ldots$ of representant
systems $S_i$ of the quotient groups $\Gamma/\Gamma_i.$
If ${\mathcal F}={\mathcal F}(M)$ is the  fundamental domain of
$\overline M\to M,$ then ${\mathcal F}(M_i):={\mathcal F}S_i$ defines fundamental domains
for all $\overline M\to M_i$.  For the classical Eta-invariant we have from Proposition \ref{s3.10}

\begin{equation*}
\begin{split}
\eta(D_i) & = \int^\infty_0 \frac{1}{\sqrt{\pi t}} \int_{{\mathcal F}(M_i)}
\sum_{\gamma\in\Gamma_i} tr_E([\overline D e^{-t\overline D^2}](x,x\gamma))dx
dt  \\
&  = d_i\int^\infty_0 \frac{1}{\sqrt{\pi t}}\int_{{\mathcal F}(M_i)}
\sum_{\gamma\in\Gamma_i} tr_E([\overline D e^{-t\overline D^2}](x,x\gamma))dx
dt,  
\end{split}
\end{equation*}

thus

\begin{equation}\label{20}
\frac{1}{d_i}\eta(D_i)-\eta_\Gamma(\overline D)=\int^\infty_0\frac{1}{\sqrt{
\pi t}}\int_{{\mathcal F}(M)}\sum_{e\neq\gamma\in\Gamma_i}
tr_E ([\overline D e^{-t\overline D^2}](x,x\gamma))dx dt.
\end{equation}

Formula (\ref{20}) raises the question of the convergence of the  $\frac{1}{d_i}\eta
(D_i)$ for $i\to\infty$. The following theorem gives a partial answer.

\begin{theo}\label{th3.12} {\bf (\cite{ChG})}Assume that one of the following conditions hold:
\begin{enumerate}
\item[(a)] $\ker \overline D=\{0\}.$
\item[(b)] $\overline D^2$ is the Laplace operator on $\Lambda T^*\overline M.
$   
\end{enumerate}

Then $\lim_{i\to\infty}\frac{1}{d_i}\eta(D_i)=\eta_\Gamma(\overline D).$
\end{theo}

{\it Proof.} We split the integral  in (\ref{20}) into integrals over
$[0,\kappa]$ and $[\kappa,\infty[.$ Then

\begin{equation*}
\begin{split}
|\frac{1}{d_i} \eta(D_i)_\kappa-\eta_\Gamma(\overline D)_\kappa| &\le
\int^\kappa_0 \frac{1}{\sqrt{\pi t}}\int_{{\mathcal F}(M)}\sum_{e\neq\gamma
\in\Gamma_i}
|tr_E ([\overline D e^{-t\overline D^2}](x,x\gamma))|dx dt  \\
&\le \int^\kappa_0 \frac{C}{\sqrt{\pi t}}\sum_{e\neq\gamma\in\Gamma_i}
e^{-(d(x_0,x_0\gamma)-r_0)^2/6c^2t}dt \\
& \le C \sum_{e\neq\gamma\in\Gamma_i}\kappa e^{-(d(x_0,x_0\gamma)-r_0)^2/6c^2
\kappa}\\
& \le C(\kappa)\sum_{e\neq\gamma\in\Gamma_i} e^{-\delta|\gamma|^2/\kappa}.
\end{split}
\end{equation*}

This series converges absolutely, such that the RHS converges to $0$ for $i\to
\infty$. Using the definition of $\eta_\Gamma(\overline D)^\kappa$ we have

$$|\eta_\Gamma(\overline D)^\kappa|\le \mbox{\rm tr}_\Gamma(e^{-\kappa\overline D^2})
-b_\Gamma  \;\mbox{and} \;
\frac{1}{d_i}\eta(D_i)^\kappa\le\frac{1}{d_i}(\mbox{\rm tr}(e^{-\kappa D^2_i})-b_i).$$

Under condition (a) Lemma \ref{l3.11} implies $\lim_{i\to\infty}b_i/d_i=b_
\Gamma.$ The same holds true under condition (b) as is shown in \cite{Lue}.
We now have

\begin{equation*}
\begin{split}
&|\frac{\eta(D)_i}{d_i}-\eta_\Gamma(\overline D)|\le C(\kappa)   \sum_{e\neq 
\gamma\in\Gamma_i}e^{-\delta|\gamma|^2/\kappa}\\
&\qquad+ \left|\frac{\mbox{\rm tr}(e^{-tD^2_i})-b_i}
{d_i}-\mbox{\rm tr}_\Gamma(e^{-t\overline D^2})+b_\Gamma\right|
  +2(\mbox{\rm tr}_\Gamma(e^{-\kappa\overline D^2})-b_\Gamma)
\end{split}
\end{equation*}

The last summand can be made arbitrarily small by choosing an appropriate  $\kappa$. 
The two remaining summands can be made arbitrarily small by letting $ i \to \infty$.

\rightline{$\blacksquare$}

\section{ Coverings of noncompact manifolds}

This Chapter contains some $\Gamma$-analogues of classical Theorems from the index theory 
on noncompact manifolds.

\subsection{Decomposition principle}

According to the classical decomposition principle (cf. \cite{DoLi}),
the essential spectrum of a differential operator on a noncompact manifold
is determined by the operator 'at infinity'.
This Section proves the analogous assertion for coverings.

For our variant of the decomposition principle, we look at hermitian vector bundles 
$E_i\to N_i$ of bounded geometry and selfadjoint operators $P_i\in \mbox{\rm UDiff}^1
(N_i,E_i).$ We assume that there is a decomposition

$$N_1=U\cup K_1,\quad N_2=U\cup K_2$$

with $K_i$ compact, and that all structures on $N_1,N_2$ are isometric over
$U$. 
We also have to assume that the pre-image $\overline U$ of $U$ 
under $\overline N_1\to
N_1$ is $\Gamma$-isometric to the pre-image of $U$ under $\overline N_2
\to N_2$.

\begin{satz}\label{s4.4} {\bf (Decomposition principle)}
$\mbox{\rm spec}_{\Gamma,e}(\overline P_1)=\mbox{\rm spec}_{\Gamma,e}(\overline P_2).$
\end{satz}

{\it Proof.} Let $\lambda\in \mbox{\rm spec}_{\Gamma,e}(\overline P_1).$ Then
for all
$\epsilon>0$ the space $G_\epsilon:=\H_{\overline P_1}(]\lambda-\epsilon,
\lambda+\epsilon[)$ has infinite $\Gamma$-dimension. We write $E(\epsilon):
=E_{\overline P_1}(]\lambda-\epsilon,\lambda+\epsilon[).$
All Sobolev-norms are equivalent on $G_\epsilon$. This follows from the estimate
for $f\in G_{\epsilon}$ and $k \in\N$

\begin{equation*}
\begin{split}
\parallel f\parallel_{H^{2k}(\overline N_1)}&  \le C  \left(
\parallel f\parallel_{L^2(\overline N_1)}
+\parallel(\overline P_1-\lambda)^k f\parallel_{L^2(\overline N_1)}\right)\\
&  \le (C+\epsilon^\kappa)\parallel f\parallel_{L^2(\overline N_1)}.
\end{split}
\end{equation*}

The proof now proceeds to show that a $\Gamma$-infinite dimensional space
of sections $f$ in $G_{\epsilon}$ essentially lives in $\overline U$.
For this, choose cut-off functions $\phi,\psi\in C^\infty_c(N_1)$ with 
$\phi|_{K_1}\equiv 1$ and $\psi|_{\supp(\phi)}\equiv 1.$ According to
\ref{s3.1} the concatenation

\begin{equation*}
C_{\overline\psi}:(G_\epsilon,\parallel\cdot\parallel_{L^2})
   \overset{Id}{\to}
(G_\epsilon,\parallel\cdot\parallel_{H^2})\overset{\overline\psi}{\to}
H^1(\overline N_1,\overline E_1)
\end{equation*}

is $\Gamma$-compact. Then the map $C^*_{\overline\psi}C_{\overline\psi} \; (\neq\overline
\psi^2)$ is $\Gamma$-compact, selfadjoint on $(G_\epsilon,\parallel\cdot
\parallel_{L^2})$, and

$$\widetilde G_\epsilon:=\H_{C^{*}_{\overline\psi}C_{\overline\psi}}
( ]-\epsilon^2,\epsilon^2[ )\subset G_\epsilon $$
is of infinite $\Gamma$-dimension.
For  $f\in\widetilde G_\epsilon$ we then find 

$$
\parallel\overline\psi f \parallel^2_{H^1(\overline N_1)} = \langle C_{\overline
\psi} f, C_{\overline\psi}f\rangle_{H^1(\overline N_1)} = 
\langle C^*_{\overline\psi}C_{\overline\psi}f,f\rangle_{L^2(\overline N_1)}
\le \epsilon^2\parallel f\parallel^2_{L^2(\overline N_1)},$$
and
\begin{equation*}
\begin{split}
& \parallel(\overline P_2-\lambda)(1-\overline\phi)f\parallel_{L^2(\overline
N_2)}  \\
& \qquad \le  
\parallel[\overline P_1,\overline\phi]f\parallel_{L^2(\overline N_1)}
+\parallel(1-\overline\phi)(\overline P_1-\lambda)f\parallel_{L^2(\overline N
_1)}   \\
& \qquad \le C\parallel\overline\psi f\parallel_{H^1(\overline N_1)}  +
\parallel(\overline P_1-\lambda)f\parallel_{L^2(\overline N_1)}\le
\epsilon(1+C)\parallel f\parallel_{L^2(\overline N_1)}.
\end{split}
\end{equation*}

This shows that $$(1-\overline\phi)\widetilde G_\epsilon
\subset \H_{\overline P_2}(]\lambda-\epsilon(1+C),\lambda+\epsilon(1+C)[).$$
But, again using Proposition \ref{s3.1}, the multiplication operator
$1-{\overline\phi}$ on $\widetilde G_\epsilon$ is $\Gamma$-Fredholm in $L^2$.  
This means that the LHS is of infinite $\Gamma$-dimension, 
i.e. $\lambda\in \mbox{\rm spec}_{\Gamma,e}(\overline P_2).$ 

\rightline{$\blacksquare$}

\section{Manifolds with cylindrical ends}

In this Chapter, $(X,g)$ is a connected (oriented, riemannian) manifold with
cylindrical ends of even dimension $n.$ Thus,  $X$is of the form

$$X=X_0\cup Z \,\mbox{with}\, Z=M\times[0,\infty[,\quad M=\partial X_0\quad
\mbox{compact}.$$
On  $X$ we have the  Clifford bundle
$(E,\nabla^E,h^E)$. We assume that all these structures are of product form over the cylinder:
$$g|_Z=g^M+\langle\cdot,\cdot\rangle_{[0,\infty[},\;
E_Z\cong \rho^*(E|_M), \,
h^E|_Z=\rho^*(h^E|_M), \,
\nabla^E|_Z= \rho^*\nabla^{E|_M}$$
etc.. Here $\rho:Z=M\times[0,\infty[\to M$ denotes the projection onto the basis of the cylinder, 
where we have identified $M$ with $M\times\{0\}\subset X$. 
We will often write $X=X_\kappa\cup Z_\kappa,$ where $X_\kappa=X_0\cup
(M\times[0,\kappa]),$ and $Z_\kappa=M\times[\kappa,\infty[$.

Choose the orientation on (all components of)  $M$ such that
$E_1,\ldots,E_{2n-1},\frac{\partial}{\partial r}$ is a positively oriented local frame in $TZ$, 
whenever $E_1\ldots,E_{2n-1}$ is a positively oriented local frame
in $TM$. Then Clifford multiplication with $\frac{\partial}
{\partial r}$ gives an isomorphism of $C(M)$-modules
$${\bold c}(\frac{\partial}{\partial r}):E^+|_M\overset{\sim}{\to}E^-|_M$$
Writing $F$ for the $C(M)$-module $E^+|_M$ we can thus identify 
$E|_Z$ and $\rho^*(F\oplus F)$ as Clifford modules via

$$
{\bold c}^E(X)\overset{\wedge}{=}
\begin{pmatrix}
0  &  {\bold c}^F(W)  \\
{\bold c}^F(W)  &  0
\end{pmatrix},\, \mbox{for}\, W  \in TM,\quad{\bold c}^E(\frac{\partial}{\partial 
r
})\overset{\wedge}{=}
\begin{pmatrix}
0  & -1  \\
1  &  0
\end{pmatrix}.$$

(Cf. Annex A). In this representation the Dirac operator $D$ on $\rho^*(F\oplus F)$ over $Z$ 
has the form

\begin{equation}\label{22}
D={\bold c}(\frac{\partial}{\partial r})\frac{\partial}{\partial r}+{\bold c}|_
M
\circ \nabla^{E|_M}={\bold c} (\frac{\partial}{\partial r})\frac{\partial}
{\partial r}+\Omega D^F,\;
\Omega=\begin{pmatrix}
0 & 1 \\
1 & 0
\end{pmatrix}.
\end{equation}

Here and in the future, we drop mention of the identifying map $\rho^*$ and just write $E|_Z\cong F\oplus F.$

In this Chapter, we consider the $\Gamma$-index theory of the lifted
Dirac operator $\overline D$ on a regular covering $\overline X$ of $X$,
with lifted structures $\overline Z, \overline E$ etc.. 
Section 6.3 contains the proof of the L$^2$-$\Gamma$-index theorem
for $\overline D$ using methods from \cite{Me}, \cite{Mue2}. 
This will require the introduction of a spectral modification of $\overline D^F$ 
which effaces the spectrum of this operator around $0$. 
The resulting modification $\overline D_{\epsilon,u}$ of $\overline D$ then is $\Gamma$-Fredholm, 
but not a differential operator. In the first two Sections of we therefore
go into some of the details of the analysis of such operators.

\subsection{$\Gamma$-operators with product structure on the cylinder}\label{subsecprodukt}

Let $T\in Op^1_\Gamma(\overline X,\overline E)$ be elliptic and selfadjoint.
 $T$ is said to have {\it product structure}, if its restriction to $\overline X_0$
is a differential operator $\overline P|_{\overline X_0}$, lifted from a
uniformly elliptic and formally selfadjoint differential operator 
$\overline P\in \mbox{\rm UDiff}^1(\overline X,\overline E)$,
and its restriction to $Z$ looks like

\begin{equation}\label{23}
T={\bold c}(\frac{\partial}{\partial r})\frac{\partial}{\partial r}+
\Omega B(r)=
\begin{pmatrix}
0 & B(r)-\frac{\partial}{\partial r} \\
B(r)+\frac{\partial}{\partial r}  & 0  
\end{pmatrix}.
\end{equation}

Here, $[0,\infty[\ni r\longmapsto B(r)\in Op^1_\Gamma(\overline M,
\overline F)$ is a  differentiable family of elliptic and selfadjoint operators. 
We assume that  $B(r)\equiv \overline B_0$ for $ r < 1$ and  $B(r)\equiv B$ 
for $r > 2$.
We recall that we allow $\overline M, \overline Z$ to have countably many components.

Let's start the analysis of these operators by noting that the results
of Section \ref{subsecprop}, especially (\ref{11}) and
(\ref{12}) remain true for $T$ on $\overline X_0$, and that 
for each $\xi_0\in C^\infty_c(\overline X,\overline E)$ 
a  {\it unique} solution $\xi(t):=e^{itT}\xi_0$ of the wave equation for $T$ exists.
The following Lemma shows that we still have an energy estimate for $\xi(t)$, though only
along the cylinder:

\begin{lem}\label{l5.1} Choose  $U=\overline M\times]a,b[\;\subset \overline Z$ with
$0<a<b.$ For $\Lambda<a$ the norm $\parallel\xi(t)\parallel_{L^2(B(U,
\Lambda-t))}$ is then  monotonously decreasing in $t.$ 
Especially,  $e^{itT}$ has propagation speed along the cylinder $\le 1$, since
 $ \supp(\xi_0)\subset U$ 
implies  $\supp(\xi(t))\subset B(U,t).$
\end{lem}

{\it Proof.} The proof follows the proof of Proposition 5.5 in \cite{Ro1}:

\begin{equation*}
\begin{split}
&\frac{\partial}{\partial t}\parallel\xi(t)\parallel^2_{L^2(B(U,\Lambda-t))}
= \frac{\partial}{\partial t}\int_{B(U,\Lambda-t)}|\xi(t)|^2(z)dz  \\
& \le |\int_{B(U,\Lambda-t)}(\langle \xi(t),iT\xi(t)\rangle
  +\langle iT\xi(t),\xi(t) \rangle)(z)dz|-\int_{\partial B(U,\Lambda-t)}
|\xi(t)|^2(z)dz
\end{split}
\end{equation*}

Now $T={\bold c}(\frac{\partial}{\partial r})\frac{\partial}{\partial r}
+\Omega B(r),$ and the domain of integration  $B(U,\Lambda-t)$ is the product
of $\overline M$ and an intervall. From the selfadjointness of $B(r)$over $M$ we can deduce

$$\int_{\overline M}(\langle i\Omega B(r)\xi(t),\xi(t)\rangle(x)+
\langle\xi(t),i\Omega B(r)\xi(t)\rangle(x))dx=0.$$

The estimation of the derivative above can therefore be continued as follows:
\begin{equation*}
\begin{split}
&\frac{\partial}{\partial t}\parallel\xi(t)\parallel^2_{L^2(B(U,\Lambda-t))}\\
&\le|\int_{B(U,\Lambda-t)}\frac{\partial}{\partial r}\langle\xi(t),{\bold c}
(\frac{\partial}{\partial r})\xi(t)\rangle(z)dz|-\int_{\partial B(U,
\Lambda-t)} |\xi(t)|^2(z)dz.
\end{split}
\end{equation*}

Since $|{\bold c}(\frac{\partial}{\partial r})|=1$ this must be less than or equal $0.$

\rightline{$\blacksquare$}

This Lemma can be used to apply the methods from Section 2.3 for the operator
family $T$ over $\overline Z$.  For simplicity, continue the
 $r\mapsto B(r)$ by the constant operator $B_0$ on $r\in \R_-$, thus obtaining a family over all of 
 ${\R}$.
Now define the reference operator ${\mathcal S}={\bold c}(\frac{\partial}
{\partial r}) \frac{\partial}{\partial r}+\Omega B(r)$ on the cylinder
$\overline F\oplus\overline F\to\overline M\times {\R}.$

\begin{satz}\label{s5.2} Choose $z_i=(x_i,s_i)\in\overline Z=\overline M\times
[0,\infty[$, $s_1, s_2>r_1,$ $r_1=r_1(\overline M\times{\R})$ as in Propositions
\ref{s2.8} and \ref{s2.9}.

\begin{enumerate}
\item[(a)] If $|s_1-s_2|>2r_1$ and $t\in{\R}_+$, then

$$|\nabla^l_{z_1}\nabla^k_{z_2}[T^m e^{-tT^2}](z_1,z_2)|\le C(k,l,m,T)
e^{-(|s_1-s_2|-r_1)^2/6t}.$$
\item[(b)] Choose cut-off functions $\psi_1,\psi_2\in C^\infty(\overline Z)^\Gamma$ 
whose supports in the $r$-direction have distance $d$.
Then for $t\in{\R}_+$ the following estimate of the operator norms holds:
$$\parallel\psi_1 T^me^{-tT^2}\psi_2\parallel\le C(m,\psi_i)e^{-d^2/6t}.$$
Note that this estimate that the constants are independent of $T$ as long as
$T$ has propagtion speed $ \leq 1$ along the cylinder.
\item[(c)] The relative version of b) also holds for $t\in{\R}_+$,  $s_1, s_2>r_1$:

\begin{equation*}
\begin{split}
|\nabla^l_{^z1}\nabla^k_{^z2}([T^me^{-tT^2}](z_1,z_2)&  -
[{\mathcal S}^me^{-t{\mathcal S}^2}](z_1,z_2))| \\
&  \le C(k,l,m,T)e^{-((s_1\wedge s_2)-r_1)^2/6t}.
\end{split}
\end{equation*}

\end{enumerate}
\end{satz}

\rightline{$\blacksquare$}

Using a variant of the Decomposition Principle \ref{s4.4},
we now describe the $\Gamma$-essential spectrum of $T^2$ on $\overline X$.
The reference operator at infinity on $\overline M \times {\R}$ is
$S={\bold c}
(\frac{\partial}{\partial r})\frac{\partial}{\partial r}+\Omega B.$ 
This will be compared with $T\sqcup-T$ on $\overline X\sqcup
 -\overline X$:

\begin{lem}\label{l5.3}$\mbox{\rm spec}_{\Gamma,e}(T^2)=\mbox{\rm spec}_{\Gamma,e}(S^2).$
\end{lem}

\rightline{$\blacksquare$}

To calculate the $\Gamma$-essential spectrum of $T$, it is therefore
enough to calculate the $\Gamma$-essential spectrum of $S^2=B^2-(\frac{\partial}{\partial r})
^2$. This can obviously be done on just one 'half' $L^2
(\overline M\times{\R},\overline F)$. 
Thus choose a spectral resolution $V$ of $B,$ i.e. $V:L^2(\overline M,\overline F)
\overset{\sim}{\to}L^2({\R}_\lambda\times {\N},\mu_B(\lambda,j))$ is a unitary equivalence
such that $Vf(B)V^{-1}$ is the multiplication operator with the function
$f$ on $L^2({\R}\times {\N},\mu_B(\lambda,j)).$ 
Write  $U:L^2({\R}_r)\to L^2({\R}_y)$ for Fourier transformation along $r$,
to obtain the unitary equivalence:
$$W=V\otimes U: L^2(\overline 
M,\overline F)\otimes L^2({\R}_r)\to L^2(({\R}_\lambda\times{\N})\times 
{\R}_y,\mu_B\times dy),$$
which transforms the action of  $S$ into a multiplication operator:
$$WSW^{-1}  \cong \lambda^2+y^2.$$

The $\Gamma$-trace on the Hilbert $\Gamma$-module $L^2(\overline M\times {\R},
\overline F)$ can then be described as the product of $tr_{\Gamma}^{\overline M}$ on
$L^2(\overline M,\overline F)$  with the usual trace $tr^\R$ on $L^2
(\R).$ Writing $\omega:=\mbox{\rm inf}(\supp(\mu_{\Gamma,B^2}))$ we thus find for 
$0\le a<b$ with $\omega<b$ and $\epsilon<(b-\omega)/2$ 

\begin{equation*}
\begin{split}
\mu_{\Gamma,S^2} & (]a,b[)=\dim_\Gamma^{\overline M}\otimes\dim^{\R}
[ \H_{S^2}(]a,b[)]  \\
& =\dim_{\Gamma}^{\overline M}\otimes \dim^{\R}W^{-1}\left[L^2(\{(\lambda,j,y)|
a<\lambda^2+y^2<b\},\mu_B(\lambda,j)\times dy)\right]  \\
& \ge \dim_\Gamma^{\overline M}V^{-1}\left[L^2(\{(\lambda,j)|a<\lambda^2<\omega +\epsilon
\},\mu_B(\lambda,j))\right]\\
&\qquad \cdot\dim^{\R}\left[L^2(\{y|0  <y^2<\epsilon\},dy)\right]
\end{split}
\end{equation*}

which equals $\infty$ since the first factor is non-zero and the second factor
is infinite. In the same manner, one shows that $\mu_{\Gamma,S^2}(]a,b[)=0$, if
$\omega\ge b.$ 

We have thus shown
\begin{satz}\label{s5.4}
$\mbox{\rm spec}_{\Gamma,e}(T^2)=[\inf(\supp(\mu_{\Gamma,B^2})),\infty[.$
\end{satz}

\rightline{$\blacksquare$}

\subsection{The $L^2$-$\Gamma$-index}\label{subsecl2index}

Proposition \ref{s5.4} implies (together with Proposition \ref{s1.19}) that
 $T\in Op^1_\Gamma (\overline X,\overline E)$, with product structure (\ref{23}), 
is $\Gamma$-Fredholm if and only if $B$ is invertible. 
In this Section, we describe the structure of the null space of such operators
without the condition of invertibility on $B$.

The analysis of the asymptotics at infinity of the sections in $\null(T^\pm)$ can be done
most naturally in the context of weighted  $L^2$- and Sobolev-spaces on $X$.
 To introduce these concepts, let $\vartheta\in C^\infty(\overline X)$ be a weight
function with 

$$\theta=\theta(r)=r\;\mbox{on}\; \overline Z_3,\,
\theta|_{\overline X_{2}}\equiv 0,\quad
\vartheta=\theta^\prime(r)\;\mbox{on}\;{\overline Z_2},\;
\vartheta|_{\overline X_{2}}\equiv 0.$$

Then, the operator $T$ is closed as an operator on $e^{u\theta}L^2,$ $u\in{\R}$ with
domain of definition $e^{u\theta}H^1$ and we have
a commutative diagram:

\begin{equation}\label{25}
\begin{split}
e^{u\theta} H^1({\overline X},\overline E^{\pm})
&  \overset{T^{\pm}}{\longrightarrow} e^{u\theta}L^2({\overline X},\overline E^{\mp})  \\
\uparrow e^{u\theta} &     \qquad\qquad\quad     \uparrow e^{u\theta}\\
 H^1({\overline X},\overline E^{\pm})
&   \overset{T^{\pm}\pm u\vartheta}{\longrightarrow}
L^2({\overline X},\overline E^{\mp})
\end{split}
\end{equation}

The  vertical maps in this diagram are isomorphisms of Hilbert $\Gamma$-modules.
The following Lemma states  the most  important properties of these weighted Sobolev spaces.
\begin{lem}\label{weighted}
 Choose $\delta^\prime<\delta \in \R$. Then
\begin{enumerate}
\item[(a)] The map $\iota:e^{-\delta\theta}L^2({\overline X})
  \hookrightarrow e^{-\delta^\prime \theta}L^2({\overline X})$ is a continuous embedding.
\item[(b)] The map $\iota:e^{-\delta\theta}H^1({\overline X})
  \hookrightarrow e^{-\delta^\prime \theta}L^2({\overline X})$ is $\Gamma$-compact.
\item[(c)] The subspace ${\mathcal W}$ is $\Gamma$ finite-dimensional in $e^{-\delta^\prime \theta}L^2$, if and only if
 $\iota({\mathcal W})$ is of finite $\Gamma$-dimension in $e^{-\delta \theta}L^2$. In that case
$$ \dim_\Gamma({\mathcal W}\subset e^{-\delta^\prime \theta}L^2)=\dim_\Gamma(\iota({\mathcal W})\subset e^{-\delta \theta}L^2) $$
\end{enumerate}
\end{lem}
{\it Proof.}
(a) is clear, for (b) set $\delta^\prime=0$ wlog. First note that the restriction operator 
$$\Lambda_\kappa:e^{-\delta \theta}H^1(\overline X)\overset{\cdot \chi_{\overline X_\kappa}}{\longrightarrow} L^2(\overline X)$$
is $\Gamma$-compact according to Rellich's Theorem \ref{s3.1}. 
We proceed to show that the operators $\Lambda_\kappa$ converge to $1$ in norm.
For $\xi \in e^{-\delta \theta}H^1$:

\begin{equation*}
\begin{split}
 &\parallel(1-\Lambda_\kappa)\xi\parallel^2_{L^2(\overline X)}  =\int^\infty_\kappa\parallel\xi(\bullet,r)\parallel^2_{L^2
(\overline M)}dr  \\
&  \le  e^{-2\delta \kappa} \int^\infty_\kappa e^{2\delta r}\parallel\xi(\bullet,r)\parallel^2_{L^2
(\overline M)}dr  
 \le  e^{-2\delta \kappa} \parallel \xi \parallel_{e^{-\delta \theta} L^2}^2 \le  e^{-2\delta \kappa} 
\parallel \xi \parallel_{e^{-\delta \theta} H^1}^2.
\end{split}
\end{equation*}
This  converges to $0$ for $\kappa \to \infty$ proving (b).

For (c) it suffices to note that
 $$\iota:({\mathcal W}\subset e^{-\delta \theta}L^2)\to (\iota({\mathcal W})\subset e^{-\delta^\prime \theta}L^2) $$
is a quasiisomorphism and to recall Lemma \ref{l1.8}.

\rightline{$\blacksquare$}

Part (c) of this Lemma enables us to indiscriminately use the notation  $\dim_\Gamma$ on all weighted
Sobolev spaces.
The Lemma also implies
\begin{cor}\label{l5.5} As before let $T\in Op^1_\Gamma
(\overline X,\overline E)$, elliptic, selfadjoint and with product structure
(\ref{23}) and let $u \in \R$. Then $e^{u\theta}L^2\mbox{\rm-null}(T)$ is of finite $\Gamma$-dimension.
\end{cor}

{\it Proof.} Again wlog set $u=0$. As all Sobolev-norms are equivalent on $\mbox{\rm null}(T)$,
Lemma \ref{weighted}(b) implies that the map

$$\iota: (\ker(T), L^2) \longrightarrow (\ker(T), H^1) \longrightarrow e^{\delta\theta }L^2$$

is $\Gamma$-compact for all $\delta > 0$. Hence
$\iota(\ker(T), L^2)\subset e^{\delta\theta }L^2$ is of finite $\Gamma$-dimension.
Part (c) of Lemma \ref{weighted}  then implies that $\mbox{\rm null}(T) \subset L^2$
 is of finite $\Gamma$-dimension.

\rightline{$\blacksquare$}

\begin{Def}\label{d5.6} Let $T\in Op^1_\Gamma(\overline X,\overline E)$ 
elliptic, selfadjoint and with product structure (23). The
$L^2\mbox{{\rm-}}\Gamma$-index of $T$ is defined as 
$$L^2\mbox{\rm-ind}_\Gamma(T):=\dim_\Gamma(\ker(T^+))
-\dim(\ker(T^-)).$$
\end{Def}

We now introduce the following modification of the  Dirac operator 
on $X$:
$$\overline D_{\epsilon,u}:=
\overline D+\vartheta\Omega(u-A\Pi_\epsilon), \quad A=\overline D^F, 
\quad\Pi_\epsilon=E_A(]-\epsilon,\epsilon[).$$
We also agree to write
$\overline D_\epsilon$ for $\overline D_{\epsilon,0}.$ 
The operator
$\overline D_{\epsilon,u}$ is our prototype of an eliptic selfadjoint operator
 in $Op^1_\Gamma(\overline X,\overline E)$ with product structure over the cylinder
that we have considered above.
Its restriction to the basis $M$ of the cylinder
is $A_{\epsilon,u}=A(1-\Pi_\epsilon)+u.$ For $0<\epsilon$,  $0$ 
is an isolated point in the spectrum of $A_\epsilon:=A_{\epsilon,0},$ and
 $A_{\epsilon,u}$ is invertible for $0<|u|<\epsilon$. In this case 
the operator $\overline D_{\epsilon,u}$ is  $\Gamma$-
Fredholm according to the results in Section \ref{subsecprodukt}.
In the remainder of this Section, we establish the relationship between the 
$\Gamma$-index of the operators $\overline D_{\epsilon,u}$ and the
L$^2$-$\Gamma$-index of $\overline D.$

On the cylinder, $\overline D_{\epsilon,u}$ can be written as

$$\overline D_{\epsilon,u}={\bold c}(\frac{\partial}{\partial r})
\frac{\partial}{\partial r}+\Omega A+\vartheta\Omega(u-A\Pi_\epsilon),$$

i.e. the Sections $\xi^\pm\in e^{\infty\theta}L^2:=\bigcup_{\delta>0}
e^{\delta\theta}L^2$ in the $C^\infty$-null space of $\overline D^{\pm}_{\epsilon,u}$
satisfy the following equation over $\overline Z$

$$\left(\pm\frac{\partial}{\partial r}+\lambda+\vartheta(r)(u-\chi_\epsilon
(\lambda)\lambda)\right)V\xi^\pm=0,\quad \chi_\epsilon:=\chi_{]-\epsilon,\epsilon[}.$$
Here, $V:L^2(\overline M , \overline F )\to L^2(\mu_A)$ is again the spectral resolution of  $A$.
We thus have
\begin{equation}\label{26}
V\xi^\pm(\lambda,r)=\zeta^\pm(\lambda,i)e^{\mp u\theta(r)} e^{\mp\lambda(r-
\theta(r)\chi_\epsilon(\lambda))}
\end{equation}

for suitably chosen $\chi^\pm(\lambda,i)\in L^2(\mu_A).$ 

All solutions of  $\overline D^\pm_{\epsilon,u}\xi=0$ are thus
exponentially decreasing, constant or increasing along the cylinder.
For $\epsilon\ge 0$ , define the space of extended $L^2$-Sections in the null space of
$\overline D_\epsilon$ and its $\Gamma$-dimension $h^\pm_{\Gamma,\epsilon}$ by

\begin{equation}\label{hgamma}
\begin{split}
\mbox{\rm Ext}(\overline D^{\pm}_{\epsilon})&:=\bigcap_{u>0}e^{u\theta}L^2\mbox{\rm-ker}(\overline
D^{\pm}_{\epsilon}),\\
h^\pm_{\Gamma,\epsilon}&:=\dim_\Gamma(\mbox{\rm Ext}(\overline D^\pm_\epsilon))
-\dim_\Gamma(L^2\mbox{\rm-ker}(\overline D^\pm_\epsilon)).
\end{split}
\end{equation}
From (\ref{26}),  (\ref{25}) one deduces for $0<u< \epsilon$

\begin{equation}\label{27}
\begin{split}
L^2\mbox{\rm-ker}(\overline D^\pm_\epsilon) &= e^{-u\theta}L^2\mbox{\rm-ker}(\overline D^\pm_\epsilon)
= L^2\mbox{\rm-ker}(\overline D^\pm_{\epsilon,\mp u})\\
\mbox{\rm Ext}(\overline D^\pm_\epsilon) &=
e^{u\theta}L^2\mbox{\rm-ker}(\overline D^\pm_\epsilon)=
L^2\mbox{\rm-ker}(\overline D^\pm_{\epsilon,\pm u})
\end{split}
\end{equation}

Taking the limit $u\to 0$ is thus harmless:
\begin{lem}\label{l5.7} Let $\epsilon>0.$ Then

\begin{enumerate}
\item[(a)] $\dim_\Gamma(\ker(\overline D^\pm_{\epsilon}))=\lim_{u\searrow 0}\dim_
\Gamma(\ker(\overline D^\pm_{\epsilon,\mp u}))$\\ 
${}\qquad\qquad\qquad\qquad\qquad =
\lim_{u\searrow 0}\dim_\Gamma(\ker(\overline D^\pm_{\epsilon,\pm u}))-h^\pm_{\Gamma,
\epsilon}.$
\item[(b)] $L^2\mbox{\rm-ind}_\Gamma(\overline D_\epsilon)=\lim_{u\searrow 0}\mbox{\rm ind}_
\Gamma(\overline D_{\epsilon,u})-h^+_{\Gamma,\epsilon}$\\ 
$ {}\qquad\qquad\qquad\qquad \qquad=
\lim_{u\searrow 0}\mbox{\rm ind}_\Gamma(\overline D_{\epsilon,-u})+h^-_{\Gamma,
\epsilon}$
\end{enumerate}
\end{lem}

{\it Proof.} These claims follow from (\ref{27}), diagram (\ref{25}) and Lemma \ref{weighted}.

\rightline{$\blacksquare$}

The description of the null spaces of $\overline D_\epsilon$
for $\epsilon \searrow 0$ is a little bit more subtle.
Let's start by collecting some of the consequences of (\ref{26}):

\begin{lem}\label{null spacemap} 
Choose $\epsilon >\delta >0$, $\delta^\prime \in \R$. Then
\begin{enumerate}
\item[(a)] $\xi \in e^{\delta^\prime \theta}L^2\mbox{\rm-ker}(\overline D^+) \Rightarrow
 \xi |_{\overline Z}= e^{-rA }\zeta,\quad \zeta \in{\mathcal H}_A(]-\delta^\prime,\infty[)$ 
\item[(b)] $\xi \in L^2\mbox{\rm-ker}(\overline D^+_\epsilon) \Rightarrow
 \xi |_{\overline Z}= e^{-rA +\theta(r)A\Pi_\epsilon}\zeta,\quad \zeta \in{\mathcal H}_A(]\epsilon,\infty[)$ 
\item[(c)] $\xi \in e^{\delta \theta}L^2\mbox{\rm-ker}(\overline D^+_\epsilon) \Rightarrow
 \xi |_{\overline Z}= e^{-rA +\theta(r)A\Pi_\epsilon}\zeta,\quad \zeta \in{\mathcal H}_A(]-\epsilon,\infty[)$ 
\item[(d)] The operator
$$e^{\pm \theta(r)A\Pi_\epsilon}:L^2(\overline X, \overline F)\rightarrow
e^{2\epsilon \theta}L^2(\overline X, \overline F)$$ 
is quasiisometric onto its image.
\item[(e)] $D^\pm e^{\mp \theta(r)A\Pi_\epsilon}=e^{\mp \theta(r)A\Pi_\epsilon} D^\pm_\epsilon$
\end{enumerate}
Analogous statements can be made for $\overline D^-$.
\end{lem}

{\it Proof.} These claims all follow from the representation 
$$ \overline D_\epsilon^\pm = A \pm \frac{\partial}{\partial r}-\vartheta(r)A\Pi_\epsilon 
\quad\mbox{on} \quad\overline Z,$$
and the description of solutions in (\ref{26}). For (d), note in addition that the
operator  $e^{\pm \theta(r)A\Pi_\epsilon}$ has no null space.

\rightline{$\blacksquare$}
 
\begin{lem}\label{l5.8} 
\begin{enumerate}
\item[(a)] 
$\lim_{\epsilon\searrow 0}\dim_\Gamma(\ker(\overline D^\pm_\epsilon))=\dim_\Gamma
(\ker(\overline D^\pm)).$
\item[(b)] $\lim_{\epsilon\searrow 0}L^2\mbox{\rm-ind}_\Gamma(\overline D_\epsilon)
=L^2\mbox{\rm-ind}_\Gamma(
\overline D).$
\item[(c)] $\lim_{\epsilon\searrow 0}\dim_\Gamma(\mbox{\rm Ext}(\overline D^\pm_\epsilon))=\dim_
\Gamma(\mbox{\rm Ext}(\overline D^\pm)).$
\end{enumerate}
\end{lem}

{\it Proof.} For (a), let $\xi \in L^2\mbox{\rm-ker}(\overline D^+_\epsilon)$, thus on the cylinder
$ \xi |_{\overline Z}= e^{-rA+\theta A\Pi_\epsilon }\zeta$ with $\zeta \in{\mathcal H}_A([\epsilon,\infty[)$. 
Here we then also have

\begin{equation*}
\begin{split} 
\overline D^+ \xi |_{\overline Z} & =  (\overline D^+_\epsilon+\vartheta(r)A\Pi_\epsilon)\xi|_{\overline Z}
= \vartheta(r)A\Pi_\epsilon(\xi|_{\overline Z}) \\
&=\vartheta(r)A\Pi_\epsilon (e^{-rA+\theta A\Pi_\epsilon }\zeta) = 0, \quad\mbox{da}\quad \Pi_\epsilon \zeta =0,
\end{split}
\end{equation*}
thus $ L^2\mbox{\rm-ker}(\overline D^+_\epsilon)\hookrightarrow  L^2\mbox{\rm-ker}(\overline D^+)$.
The operator $\overline D^+_\epsilon$, when restricted to $ L^2\mbox{\rm-ker}(\overline D^+)$,
has null space $L^2\mbox{\rm-ker}(\overline D^+_\epsilon)$ and its image satisfies 
$$\overline D^+_\epsilon (L^2\mbox{\rm-null}(\overline D^+)) =-\vartheta(r)A\Pi_\epsilon  (L^2\mbox{\rm-ker}(\overline D^+))
\subset -\vartheta(r)A e^{-rA} {\mathcal H}_A(]0,\epsilon[),$$
hence $\dim_\Gamma(\overline D^+_\epsilon (L^2\mbox{\rm-ker}(\overline D^+)))\to 0 $ for $\epsilon \searrow 0$.
This and Lemma \ref{l1.8} prove the result.

(c): Let $0 < \delta < \epsilon$. Using the description of Sections in the null space over the cylinder in \ref{null spacemap}, 
one shows that the operator $\Psi^\pm_\epsilon= e^{\pm \theta(r)A\Pi_\epsilon}$ satisfies:
$$ e^{\theta \delta}L^2\mbox{\rm-ker}(\overline D^+) \overset{\Psi^+_\epsilon}{\longrightarrow}
 e^{\theta \delta}L^2\mbox{\rm-ker}(\overline D^+_\epsilon) \overset{\Psi^-_\epsilon}{\longrightarrow}
  e^{\theta (\delta+\epsilon)}L^2\mbox{\rm-ker}(\overline D^+).$$

These maps are injective, and when restricted (e.g.) to
$$ e^{\theta \delta}L^2(\overline X,\overline F)\overset{\Psi^+_\epsilon}{\longrightarrow}
 e^{2\theta \epsilon}L^2(\overline X,\overline F) \overset{\Psi^-_\epsilon}{\longrightarrow}
 e^{4\theta \epsilon}L^2(\overline X,\overline F)
$$
also continuous.Taking $\epsilon \searrow 0$ and using Lemma \ref{weighted} then shows the result.

\rightline{$\blacksquare$}

\subsection{The L$^2$-$\Gamma$- index theorem}

In this Section we show (cf. \cite{BGV}, Chapter 3 for notation):

\begin{theo}\label{th5.11} {\bf (L$^2$-$\Gamma$-index theorem) }
$$L^2\mbox{\rm-ind}_\Gamma(\overline D)=\int_X\widehat A(X)
Ch(E/S)+\frac{1}{2}(\eta_\Gamma(A)-h^-_\Gamma+h^+_\Gamma)$$
\end{theo}

Using  Lemma \ref{l5.7} and Lemma \ref{l5.8} we will reduce the calculation of the $L^2$-$\Gamma$-index
of $\overline D$ to the calculation of the $\Gamma$-
index of $\overline D_{\epsilon,u},$ for small $0<|u|<\epsilon$.

According to Proposition \ref{s2.10} for $t \to \infty$, the operator $[e^{-t\overline D^2_{\epsilon, u}}]$
converges in $C^\infty$ to  the Schwartz-kernel $[N(\overline D_{\epsilon,u})]$ of the projection 
onto the null space of $\overline D_{\epsilon,u}.$ 
We will use cut-off functions $\phi_\kappa\in C^\infty(\overline X)^
\Gamma$ with $\phi_\kappa|\overline X_\kappa=1$ and $\phi_\kappa|_{\overline
Z_{\kappa+1}}=0,$, such that the operators $\phi_\kappa e^{-s\overline
D^2_{\epsilon,u}}\phi_\kappa$ are $\Gamma$-trace class. The $\Gamma$-index 
of $\overline D_{\epsilon,u}$ can then be calculated as

\begin{equation}\label{29}
\begin{split}
& \mbox{\rm ind}_\Gamma(\overline D_{\epsilon,u}) 
= \mbox{\rm str}_\Gamma(N(\overline D_{\epsilon,u}))
= \lim_{\kappa\to\infty}\lim_{t\to\infty} \mbox{\rm str}_\Gamma(\phi_\kappa e^{-t\overline 
D^2_{\epsilon,u}}\phi_\kappa)  \\
& = \lim_{\kappa\to\infty} \left(  \mbox{\rm str}_\Gamma(\phi_\kappa e^{-s\overline 
D^{2}_{\epsilon,u}} \phi_\kappa)-\int^\infty_s  \mbox{\rm str}_\Gamma(\phi_\kappa \overline
D^{2}_{\epsilon,u} e^{-t\overline D^{2}_{\epsilon,u}}  \phi_\kappa)dt\right).
\end{split}
\end{equation}

Of course, the RHS is independent of $s>0.$ The integral can be split as follows 

\begin{equation}\label{30}
\int^{{\kappa}}_s   \mbox{\rm str}_\Gamma(\phi_\kappa \overline
D^2_{\epsilon,u}e^{-t\overline D^2_{\epsilon,u}}\phi_\kappa)dt+
\int^{\infty}_{{\kappa}}   \mbox{\rm str}_\Gamma(\phi_\kappa \overline
D^2_{\epsilon,u}e^{-t\overline D^2_{\epsilon,u}}\phi_\kappa)dt.
\end{equation}

We first show that the second  integral in \ref{30} vanishes for $\kappa\to\infty$: Since
 $\overline D_{\epsilon,u}$ is $\Gamma$-Fredholm ist, there is  $\alpha=\alpha(u)>0$ (e.g. $\alpha
(u)=|u|/2),$ such that the projection $H_\alpha=E_{\overline D_{\epsilon,u}}
([-\alpha,\alpha])$ has finite  $\Gamma$-trace.
This implies (Cf. Remark \ref{bem1.11}:

\begin{equation*}
\begin{split}
|\int^\infty_{{\kappa}}&   \mbox{\rm str}_\Gamma(\phi_\kappa \overline
D^2_{\epsilon,u}e^{-t\overline D^2_{\epsilon,u}}\phi_\kappa)dt|  \\
& \le \int^\infty_{{\kappa}}| \mbox{\rm str}_\Gamma\left(  \phi_\kappa \overline
D^2_{\epsilon,u}e^{-\overline D^2_{\epsilon,u}/2}
(1-H_\alpha)  e^{-(t-1)\overline D^2_{\epsilon,u}}
e^{-\overline D^{2}_{\epsilon,u}/2}
\overline D_{\epsilon,u}\phi_\kappa\right)|dt   \\
& \qquad+\int^\infty_{{\kappa}}| \mbox{\rm str}_\Gamma\left(
e^{-t\overline D^2_{\epsilon,u}/2} H_\alpha\overline D_{\epsilon,u}
\phi^2_\kappa\overline D_{\epsilon,u}H_\alpha 
e^{-t\overline D^{2}_{\epsilon,u}/2}\right) |dt.  \\
& \le \int^\infty_{{\kappa}}e^{-(t-1)\alpha}| \mbox{\rm str}_\Gamma\left(
\phi_\kappa \overline D^{2}_{\epsilon,u}e^{-\overline D^2_{\epsilon,u}}
\phi_\kappa \right)|dt \\
&  \qquad+\int^\infty_{{\kappa}}| \mbox{\rm str}_\Gamma 
  \left( \overline D^{2}_{\epsilon,u}
e^{-t\overline D^{2}_{\epsilon,u}}H_\alpha\right)|dt
\end{split}
\end{equation*}

The Schwartz-kernel of $\overline D^{2}_{\epsilon,u}e^{-\overline D^{2}_{\epsilon,u}}$ 
is $UC^\infty.$ Thus the trace in the first integral is majored by $c_1+c_2\kappa,$ 
and the integral converges to $0$ for 
$\kappa\to\infty$.
Using the $\Gamma$-Fredholm property of
$\overline D_{\epsilon,u}$ we get for the second integral

\begin{equation*}
\begin{split}
\int^\infty_{{\kappa}} & | \mbox{\rm str}_\Gamma \left( \overline D^{2}_{\epsilon,
u} 
e^{-t\overline D^{2}_{\epsilon,u}} H_\alpha\right)|dt  \le 
\int^{\infty}_{{\kappa}}\int^{\alpha}_{-\alpha} x^2 e^{-tx^2}d\mu_{\Gamma,
\overline D_{\epsilon,u}}(x)dt  \\
& = \int^\alpha_{-\alpha}  e^{-{\kappa}x^2} \int^\infty_0
x^2 e^{-tx^2} dt\,d\mu_{\Gamma,\overline D_{\epsilon,u}}(x)\le C
\int^\alpha_{-\alpha} e^{-{\kappa}x^2} d\mu_{\Gamma,\overline D_
{\epsilon,u}}(x)  \\
& \le C\mu_{\Gamma,\overline D_{\epsilon,u}}([-\alpha,\alpha]).
\end{split}
\end{equation*}

Thus, the penultimate term also converges to $0$ for $\kappa\to\infty$  from Levi's theorem.

We now turn our attention to the first integral in in (\ref{30}) zu. 
The integrand can be written as follows

\begin{equation}\label{31}
\begin{split}
\mbox{\rm str}_\Gamma   & (\phi_\kappa \overline D^{2}_{\epsilon,u}
e^{-t\overline D^{2}_{\epsilon,u}}  \phi_\kappa) 
= \frac{1}{2} \mbox{\rm str}_\Gamma (\phi_\kappa
[\overline D_{\epsilon,u},\overline D_{\epsilon,u}
e^{-t\overline D^{2}_{\epsilon,u}}]
\phi_\kappa)  \\
&=  \frac{1}{2} \mbox{\rm str}_\Gamma
([\overline D_{\epsilon,u},\phi_\kappa\overline D_{\epsilon,u}
e^{-t\overline D^2_{\epsilon,u}} \phi_\kappa]
-[\overline D_{\epsilon,u},\phi^2_\kappa]
\overline D_{\epsilon,u}
e^{-t\overline D^{2}_{\epsilon,u}})\\
&= \frac{1}{2} \mbox{\rm str}_\Gamma
(-[\overline D_{\epsilon,u},\phi^2_\kappa]
\overline D_{\epsilon,u}
e^{-t\overline D^{2}_{\epsilon,u}})\\
&  =-\frac{1}{2} \mbox{\rm str}_\Gamma
({\bold c}(\frac{\partial}{\partial r})\frac{\partial(\phi^2_\kappa)}{\partial
r}\overline D_{\epsilon,u}
e^{-t\overline D^{2}_{\epsilon,u}}).
\end{split}
\end{equation}

The derivatives of $\phi_\kappa$ are supported  in $\overline M\times[\kappa,
\kappa+1]$. Thus expression (\ref{31}) can be compared with the
operator $S_{\epsilon,u}={\bold c}(\frac{\partial}{\partial r})
\frac{\partial}{\partial r}+(A(1-\Pi_\epsilon)+u)\Omega$ on the cylinder
$\overline Y=\overline M\times{\R}$. Using Proposition \ref{s5.2} (c)
we get for $(x,r)\in\overline M\times [3,\infty[ \cong\overline Z_3$

$$\left( [\overline D_{\epsilon,u}e^{-t\overline D^{2}_{\epsilon,u}}]-
[S_{\epsilon,u}e^{-tS^2_{\epsilon,u}}]\right)
(x,r;x,r)\le Ce^{-(r-3)^2/6t},$$
 
i.e. both Schwartz-kernels have the same $t \to 0$-asymptotics for large $r$.
Hence

\begin{equation*}
\begin{split}
&  \int^{{\kappa}}_s|\mbox{\rm str}_\Gamma
({\bold c}(\frac{\partial}{\partial r})\frac{\partial(\phi^2_\kappa)}{\partial
r}
\overline D_{\epsilon,u}e^{-t\overline D^{2}_{\epsilon,u}})-
\mbox{\rm str}_\Gamma
({\bold c}(\frac{\partial}{\partial r})\frac{\partial(\phi^2_\kappa)}{\partial
r}
S_{\epsilon, u}e^{-tS^2_{\epsilon,u}})|dt  \\
& = \int^{{\kappa}}_s | \int_{M\times[\kappa,\kappa+1]}\pi_* \mbox{\rm str}_E
\left({\bold c}(\frac{\partial}{\partial r})\frac{\partial(\phi^2_\kappa)}
{\partial r}
[\overline D_{\epsilon,u}e^{-t\overline D^{2}_{\epsilon,u}}-
S_{\epsilon,u}e^{-tS^2_{\epsilon,u}}]\right)dxdr|dt  \\
&\le C \int^{{\kappa}}_s\int^{\kappa+1}_\kappa
e^{-(r-3)^2/6t}drdt\le C \int^{{\kappa}}_s
e^{-(\kappa-3)^2/6t}dt\\
&\le C
\int^{1/s}_{1/{\kappa}}a^{-2}
e^{-(\kappa-3)^2a/6}da  
 \le C(\kappa^2e^{-\kappa/c_1}+e^{-c_2/{s}})
\end{split}
\end{equation*}

for large $\kappa$ and small  $s.$ Thus in the first integral in (\ref{30}) ,
when looking at the asymptotics $s\to0,$ and the limit $\kappa\to\infty$, 
the operator  $\overline D_{\epsilon,u}$  can be replaced by $S_{\epsilon,u}$.

But the  integral for $S_{\epsilon,u}$ on the cylinder can be calculated. The Schwartz-kernel
of $S_{\epsilon,u}e^{-tS^2_{\epsilon,u}}$, restricted to the diagonal in
$\overline Y\times\overline Y$, is of the form

\begin{equation*}
\begin{split}
[S_{\epsilon,u}
e^{-tS^{2}_{\epsilon,u}}]
(x,r) &  =(A_{\epsilon,u}\Omega+
\mbox{${\bold c}(\frac{\partial}{\partial r})\frac{\partial}{\partial r}$})\left(
[e^{-t(A_{\epsilon,u}\Omega)^2}]
(x,y)
\frac{e^{-\frac{(r-s)^2}{4t}}}{\sqrt{4\pi t}}\right)
|_{\overset{y=x}{s=r}} \\
& = \frac{1}{\sqrt{4\pi t}}
\Omega[A_{\epsilon,u}
e^{-tA^2_{\epsilon,u}}](x,x).
\end{split}
\end{equation*}

Now (Cf. Appendix A) $str^E({\bold c}(\frac{\partial}{\partial r})\Omega
\bullet)=-2tr^F(\bullet),$ and thus 
\begin{equation*}
\begin{split}
&\int^{{\kappa}}_s \mbox{\rm str}_\Gamma
({\bold c}(\frac{\partial}{\partial r})\frac{\partial(\phi^2_\kappa)}{\partial
r}
S_{\epsilon, u}e^{-tS^2_{\epsilon,u}}) dt\\
&  =\int^{\kappa+1}_\kappa\frac{\partial(\phi^2_\kappa)}
{\partial r}dr 
\int^{{\kappa}}_s
\frac{1}{\sqrt{4\pi t}}\int_{{\mathcal F}(M)}
-2 \,tr^F
([A_{\epsilon,u}e^{-tA^2_{\epsilon,u}}](x,x))dx\,dt,\\
&  =
\int^{{\kappa}}_s
\frac{1}{\sqrt{\pi t}}\int_{{\mathcal F}(M)}
\,tr^F
([A_{\epsilon,u}e^{-tA^2_{\epsilon,u}}](x,x))dx\,dt
\end{split}
\end{equation*}

Taking the limit $\kappa \to \infty$ and selecting the constant term $\mbox{LIM}$
in the asymptotic development for $s \to 0$, this just gives the $\Gamma$-Eta-invariant $\eta_\Gamma(A_{\epsilon,u})$
 as described in Section \ref{subseceta}.
Using (\ref{29}), (\ref{31}) our present knowledge can be summarised like this

\begin{equation}\label{32}
\mbox{\rm ind}_\Gamma(\overline D_{\epsilon,u})=\lim_{\kappa\to\infty}\mbox{\rm LIM}_{s\to0}\,
\mbox{\rm str}_\Gamma\left(\phi_\kappa e^{-s\overline D^2_{\epsilon,u}}\phi_\kappa\right)
+\frac{1}{2}\eta_\Gamma(A_{\epsilon,u}), 
\end{equation}

especially the limit in the first term exists.\par
\medskip
It remains to analyse the asymptotics of the local trace $\mbox{\rm str}^E([e^{-s
\overline D^2_{\epsilon,u}}])$ $(x,x)$ for $s\to0$. For $x\in \overline X_1$ 
the asymptotic development of $\mbox{\rm str}^E([e^{-sD^2}])$ on $X$ gives
$$
LIM_{s\to0}\,\mbox{\rm str}^E([e^{-s\overline D^2_{\epsilon,u }}])(x,x)\mbox{dvol}_
{\overline X}=
\widehat A(X)Ch(E/S)(x),\quad x\in\overline X_1.$$

`Far out` on the cylinder, for instance for $(y,r)\in\overline Z$ with
$r>4$, Proposition \ref{s5.2} allows us to look at the kernel $[e^{-sS^2_{\epsilon,u}}](y,r)$
instead of $[e^{-s\overline D^2_{\epsilon,u}}](y,r)$. Since $A_{\epsilon,u}$
is invertible, we have on $\overline M$ (Cf.  Lemma \ref{l3.6})

$$0=\mbox{\rm ind}_\Gamma(A_{\epsilon,u})=\mbox{\rm str}_\Gamma(e^{-sA^2_{\epsilon,u}}),$$

independent of $s$. Thus for $4<a<b$

\begin{equation*}
\begin{split}
&  \underset{s\to0}{\mbox{\rm LIM}} \int_{\overline M\times[a,b]}\mbox{\rm str}^E([e^{-s\overline 
D^2_{\epsilon,u }}](y,r))dy\,dr \\
&=     \underset{s\to0}{\mbox{\rm LIM}}\int_{\overline M\times[a,b]}\mbox{\rm str}^E([e^{-s\overline 
S^2_{\epsilon,u }}](y,r))dy\,dr
  =    \underset{s\to0}{\mbox{\rm LIM}}\,
\frac{b-a}{\sqrt{4\pi s}}\mbox{\rm str}_\Gamma (e^{-s\overline A^2_{\epsilon,u}})=0. 
\end{split}
\end{equation*}

The part $\overline M\times[1,4]$ is a little bit more complicated, as
 $\vartheta$ is not constant here. Using Proposition \ref{s5.2} we can replace $\overline D_{\epsilon,u}$ 
by the operator $S_{\epsilon,u}:={\bold c}(\frac{\partial}{\partial r})
\frac{\partial}{\partial r}+\Omega(A+\vartheta(u-\Pi_\epsilon A))$ in that area. 
For a cut-off function $\psi_1\in C^\infty(\overline Y)^\Gamma$ with
$\psi_1|_{\overline M\times [1,4]}\equiv 1$ and support in $\overline M\times
[0,5],$ we already know from the above and (\ref{32}) that
LIM$_{s\to 0}\mbox{\rm str}_\Gamma(\psi_1 e^{-sS_{\epsilon,u}^2}\psi_1)$ exists.
We now show

\begin{lem}\label{l5.12}
$\lim_{s\to 0}\mbox{\rm str}_\Gamma(\psi_1(e^{-sS_{\epsilon,u}^2}-e^{-sS^2_u})\psi_1)=0,\quad
S_u:=S_{0,u}.$
\end{lem}

{\it Proof.} First,
$$S_{\epsilon,u}^2-S^2_u= (\vartheta \Omega \Pi_\epsilon A)^2-{\bold c}(\frac{\partial}{\partial r})
\Omega \vartheta^\prime \Pi_\epsilon A- 2\vartheta \Omega \Pi_\epsilon A S_u$$

This is an operator of finite $\Gamma$-rank in the $\overline M$-direction.
A little care is needed as it has a first dervative acting in the $\R_r$-direction.
We apply the Duhamel-method (Proposition \ref{s2.12}):

\begin{equation}\label{33}
\begin{split}
|\mbox{\rm str}_\Gamma  & (\psi_1(e^{-sS^2_u}-e^{-sS_{\epsilon,u}^2})\psi_1)|=|\mbox{\rm str}_\Gamma
(\psi_1 e^{-\delta S^2_u}  e^{-(s-\delta)S_{\epsilon,u}^2}\psi_1)|^s
_0|  \\
&=|\int^s_0 \mbox{\rm str}_\Gamma(\psi_1^2 \Pi_\epsilon e^{-\delta S^2_u}
(S_{\epsilon,u}^2-S^2_u)\Pi_\epsilon e^{-(s-\delta)S_{\epsilon,u}^2})d\delta|  \\
&\le C \int^s_0 |\mbox{\rm tr}_\Gamma( \psi_1 e^{-\delta S^2_u}\Pi_\epsilon\psi_1)|
\parallel(S_{\epsilon,u}^2-S^2_u)\Pi_\epsilon e^{-(s-\delta)S_{\epsilon,u}^2}  \parallel d\delta 
\end{split}
\end{equation}

Using the results of Section \ref{subsecprop} one finds
\begin{equation}\label{34}
\begin{split}
 |\mbox{\rm tr}_\Gamma( \psi_1 e^{-\delta S^2_u}\Pi_\epsilon\psi_1)| &\le C \; \delta^{-1/2}\\
\parallel(S_{\epsilon,u}^2-S^2_u)\Pi_\epsilon e^{-(s-\delta)S_{\epsilon,u}^2}\parallel  &\le C \; (s-\delta)^{-1/2}
\end{split}
\end{equation}

where the constants are independent of any small $|u|<\epsilon$.
Thus (\ref{33})  can be majored by
\begin{equation}\label{35}
\begin{split}
& C \int^{s/2}_0 |\mbox{\rm tr}_\Gamma( \psi_1 e^{-\delta S^2_u}\Pi_\epsilon\psi_1)| d\delta 
+ C \int^s_{s/2}
\parallel(S_{\epsilon,u}^2-S^2_u)\Pi_\epsilon e^{-(s-\delta)S_{\epsilon,u}^2}  \parallel d\delta \\
&\quad \le  C \int^{s/2}_0  \delta^{-1/2}d\delta +C \int_{s/2}^s  (s-\delta)^{-1/2}d\delta \le C \;  s^{1/2}  
\end{split}
\end{equation} 
But this converges to $0$ for  $s\to 0$.

\rightline{$\blacksquare$}

We thus have brought back the asymptotic development of the $\Gamma$-trace of the heat kernel of
 $\overline D_{\epsilon,u}$
over the critical area $\overline M \times [1,4]$ to the asymptotic development

$$\mbox{\rm str}^E([e^{-sS^2_u}])(z,z)\sim \sum_{j\in {\N}}a_j(S_u)(z)t^{(j-n)/2},$$

with  coefficients $a_j(S_u)$  differentiable in $u$. Since $S_0$ is the Dirac operator, we have
$a_j(S_0)=0$ for $j\le n/2$.
 Using Lemma \ref{l5.12} and (\ref{32}), this implies

\begin{satz}\label{s5.13} 
The $\Gamma$-index dof the $\Gamma$-Fredholm operator $\overline D_{\epsilon,u}$
is  given by
$$\mbox{\rm ind}_\Gamma(\overline D_{\epsilon,u})=\int_X
\widehat A(X)Ch(E/S)+\frac{1}{2}\eta_\Gamma(A_{\epsilon,u})+g(u),\,
\mbox{\rm with}\quad
\lim_{u\to0}g(u)=0.$$
\end{satz}
\rightline{$\blacksquare$}

Now, from Lemma \ref{l3.8} the $\Gamma$- Eta-invariant with $0<|u|<\epsilon$ satisfies
$$\eta_\Gamma(A_\epsilon)=\frac{1}{2}(\eta_\Gamma(A_{\epsilon,u})+
\eta_\Gamma(A_{\epsilon,-u})),$$
which with Lemma \ref{l5.7} and the definition (\ref{hgamma}) of $h^\pm_{\Gamma,\epsilon}$ implies

\begin{equation*}
\begin{split}
L^2\mbox{\rm-ind}_\Gamma(\overline D_\epsilon) &  =\lim_{u\searrow0}\frac{1}{2}
\left[\mbox{\rm ind}_\Gamma(\overline D_{\epsilon,u})+\mbox{\rm ind}_\Gamma(\overline D_{\epsilon,-u})
+h^-_{\Gamma,\epsilon}-h^+_{\Gamma,\epsilon}\right]  \\
& =\int_X\widehat A(X)Ch(E/S)+\frac{1}{2}(\eta_\Gamma(A_\epsilon)+h^-_
{\Gamma,\epsilon}-h^+_{\Gamma,\epsilon}). 
\end{split}
\end{equation*}
Theorem \ref{th5.11} follows from this and the observations (Cf. Lemma \ref{l5.8})

\begin{enumerate}
\item[$\bullet$] $|\eta_\Gamma(A_\epsilon)-\eta_\Gamma(A)|\le \mbox{\rm tr}_\Gamma(E_A(]-\epsilon,\epsilon[-\{0\}))\to 0$ 
for $\epsilon\searrow 0$. 
\item[$\bullet$] $h^\pm_{\Gamma,\epsilon}\to h^\pm_{\Gamma}$  for $\epsilon\searrow 0$. 
\item[$\bullet$] $L^2\mbox{\rm-ind}_\Gamma(\overline D_\epsilon)\to L^2\mbox{\rm-ind}_\Gamma(\overline D)$  for $\epsilon\searrow 0$. 

\end{enumerate}

\rightline{$\blacksquare$}

As a first simple application of Theorems \ref{th3.12} and \ref{th5.11} consider
a residually finite covering 

$$\overline X\cdots\to X_{i+1}\to X_i\cdots \to X$$

of $X$ with lifted bundles $E_i,$ Dirac operators $D_i$ etc., such that
$\overline D^F$ satisfies conditions (a) or (b) in Theorem \ref{th3.12}:

\begin{satz}\label{s5.14}{\bf (Convergence of the modified $L^2$-Index)}
$$\lim_{i\to\infty} d^{-1}_i(L^2\mbox{\rm-ind}(D_i)-
\frac{1}{2}(h^+_i-h^-_i))=L^2\mbox{\rm-ind}_\Gamma(\overline D)-\frac{1}{2}
(h^+_\Gamma-h^-_\Gamma).$$
\end{satz}

\rightline{$\blacksquare$}

\subsection{The signature operator}

In this Section we specialise the L$^2$-$\Gamma$-index Theorem \ref{th5.11} to the
case of the signature operator. The underlying Clifford bundle is $E\otimes
W=\Lambda T^*X\otimes W,$ with $W$ a flat bundle with product structure on the cylinder of $X.$ 
Clifford multiplication is given by ${\bold c}(v)=\epsilon(v)-\iota(v),$
for $v\in T^*X.$ The grading operator $\tau\equiv \tau_X=
i^{n/2}\ast_X (-1)^{|\cdot|(|\cdot|-1)/2}$ gives a  ${\Z}_2$-grading on $E$ 
and the Dirac operator $D=S_X=d_X+d^*_X$ is also called the signature operator
in this case.

We have $S^2_{\overline X}=\Delta_{\overline X}$ and from Hodge's theorem for $\overline X$, the map

$$\H^*(\overline X;\overline W):=L^2\mbox{\rm-ker}(\Delta_{\overline X})=L^2\mbox{\rm-ker}(S_
{\overline X})\to H^*_{(2)}(\overline X;\overline W),$$

from the null space of  $\Delta_{\overline X}$ into the reduced $L^2$-Cohomology of
$\overline X$ with coefficients $\overline W,$ is an isomorphism. 

For
$\dim(X)=n=4k$ we have $\ast^2_X=1$ and the $\Gamma$-signature
$\sigma_\Gamma(\overline X;\overline W)$ of $\overline X$ with coefficients
 $\overline W$ is defined as the difference of the $\Gamma$-dimensions
of the  $+1$- and $-1$-Eigenspaces of the  quadratic form $\alpha\mapsto
(\alpha,\ast_{\overline X} \alpha)$ on $\H^{2k}(\overline X;\overline W).$ The 
$\Gamma$-signature is then just the $L^2$-$\Gamma$-index of $S_{\overline X}.$

Using the identification  $\overline E=\Lambda T^*\overline M\oplus
\Lambda T^*\overline M$ via 
\begin{equation}\label{35}
\overline F:=\Lambda T^*\overline M\overset{1+\tau_X}{\to}(1+\tau_X)
\Lambda T^*\overline X=\overline E^+
\end{equation}

over the cylinder $\overline Z$,  $S_{\overline X}$ can be written

$$S_{\overline X} ={\bold c}(\frac{\partial}{\partial r})\frac{\partial}
{\partial r} + \Omega(\ast_{\overline M}d_{\overline M}-d_{\overline M}
\ast_{\overline M})(-1)^{|\cdot|(|\cdot|-i)|\cdot|/2}\equiv {\bold c}(\frac{\partial}{\partial r})\frac{\partial}
{\partial r}+\Omega A.$$

In this Section, we prove 

\begin{satz}\label{s5.15} $h^+_{\Gamma,\epsilon}
(S_{\overline X})=h^-_{\Gamma,\epsilon}(S_{\overline X})$  for all $\epsilon\ge 0.$
\end{satz}

This immediately implies
\begin{theo}\label{th5.16} {\bf ($\Gamma$-Signature theorem)} 
$$\sigma_\Gamma
(\overline X;\overline W)=rk(W)\int_XL(X)+\frac{1}{2}\eta_\Gamma(A).$$
\end{theo}

{\it Proof.} Apply Theorem \ref{th5.11} and Proposition \ref{s5.15}, recalling that 
$\widehat A(X)Ch(\Lambda T^*X\otimes W/S)=rk(W)L(X),$ ( Cf. \cite{BGV})

\rightline{$\blacksquare$}

For residually finite coverings this then implies
\begin{cor}\label{c5.17} Let the conditions of Proposition \ref{s5.14} hold.

\begin{enumerate}
\item[(a)] $\lim_{i\to\infty} d^{-1}_i\sigma(X_i;W_i)=\sigma_\Gamma(
\overline X;\overline W).$
\item[(b)] If the universal covering $\widetilde X$ of $X$ is residually finite
then $\sigma_\Gamma(\widetilde X)$ is a proper
homotopy-invariant of $X.$
\end{enumerate}
\end{cor}

{\it Proof.}  (a) follows from \ref{th5.16} and Proposition
\ref{s5.14}, (b) then follows from the homotopy invariance of the signatures $\sigma(X_i).$

\rightline{$\blacksquare$}

The proof of Proposition \ref{s5.15} is an adaptation to the $\Gamma$-case of the proof
for the classical case given in \cite{Me}. 
Thus, consider the operator $S_{\overline
X,\epsilon}$ for $\epsilon>0.$ If, instead of the identification (\ref{35}) given above,
one uses the identification $\Lambda T^*\overline X=\Lambda 
T^*\overline M \oplus \Lambda T^*\overline M\wedge dr,$ for the Clifford bundle over $\overline Z$,
then

\begin{equation}\label{2ndid}
S_{\overline X,\epsilon}=\mbox{c}(\frac{\partial}{\partial r})\frac
{\partial}{\partial r}+S_{\overline M}\oplus S_{\overline M}-\vartheta
\widehat\Pi_\epsilon S_{\overline M}\oplus S_{\overline M}
\end{equation}

where $S_{\overline M}=d_{\overline M}+d^*_{\overline M}$ and $\widehat\Pi_
\epsilon:=E_{S_{\overline M}\oplus S_{\overline M}}(]-\epsilon,\epsilon[).$
This modification of  $S_{\overline X}$ is indeed the same as the one introduced
in Section 6.2, as the intervall $]-\epsilon,\epsilon[$ cut out of the spectrum
is symmetric around $0$. For simplicity we forget about the coefficient bundle $W$.

With respect to the identification (\ref{35}) Sections $\xi^\pm\in \mbox{\rm Ext}(S^\pm_{\overline
X,\epsilon})$ can be written  on $\overline Z_3:$

$$V\xi^\pm(\lambda,r)=\zeta^\pm(\lambda,i)(\chi_{]-\epsilon,\epsilon[}
(\lambda)+(1-\chi_{]-\epsilon,\epsilon[}(\lambda))e^{\mp\lambda r})
$$ 

with suitable $\zeta^\pm(\lambda,i)\in L^2((\pm[0,\infty[)\times {\N},
\mu_A).$ The coefficient of the component of  $\xi^\pm$ that is constant in $r$ is
$\zeta^\pm(\lambda,i)\chi_{]-\epsilon,\epsilon[}(\lambda)$ and

$$V^{-1}(\zeta^\pm(\lambda,i)(\chi_{]-\epsilon,\epsilon[}(\lambda))\in
{{\mathcal H}}_{S_{\overline M}\oplus S_{\overline M}}(]-\epsilon,\epsilon[)^
\pm,$$

Here, we use that $\tau_{\overline X}$ anticommutes with $S_{\overline M}
\oplus S_{\overline M}$ and thus induces a grading on
${{\mathcal H}}_{S_{\overline M}\oplus S_{\overline M}}(]-\epsilon,\epsilon[)$. 

Using the decomposition $\Lambda T^*\overline X\cong \Lambda T^*\overline
M\oplus \Lambda T^*\overline M\wedge dr$ into forms with or without a $dr$-component,
we find
$$
{{\mathcal H}}_{S_{\overline M}\oplus S_{\overline M}}(]-\epsilon,\epsilon[)
={{\mathcal H}}_{S_{\overline M}}(]-\epsilon,\epsilon[)
\oplus {{\mathcal H}}_{S_{\overline M}}(]-\epsilon,\epsilon[),$$

and the operation of $\tau_{\overline X}$ on the RHS is just 
$\begin{pmatrix}
0       &     \tau_{\overline M}(-1)^{n-1-|\cdot |}  \\
 \tau_{\overline M}(-1)^{|\cdot |}  &  0
\end{pmatrix}$. 

 Concatenation of these isomorphisms gives

\begin{equation*}
\begin{split}
\mbox{\rm Ext}(S^\pm_{\overline X,\epsilon})&\overset{{\mathcal J}^\pm}{\longrightarrow}  
\left[{{\mathcal H}}_{S_{\overline M}}(]-\epsilon,\epsilon[)\oplus
{{\mathcal H}}_{S_{\overline M}}(]-\epsilon,\epsilon[)\right]^\pm \\
 \xi^\pm &
\longmapsto V^{-1}(\xi^\pm(\lambda,i)\chi_{]-\epsilon,\epsilon[}
(\lambda)).
\end{split}
\end{equation*}

Following  (\ref{27}), $\ker(S_{\overline X,\epsilon})$
is closed in $e^{-\delta\theta}L^2,$ for  $0<\delta<\epsilon$. Also, Ext$(S_{\overline X,\epsilon})$
is closed  in $e^{\delta\theta}L^2.$ It is than straightforward to verify that

\begin{lem}\label{l5.18} 
\begin{enumerate}
\item[(a)] The sequence
\begin{equation*}
\begin{split}
(L^2\mbox{\rm-ker} (S^\pm_{\overline X,\epsilon}),\parallel\cdot\parallel_
{e^{-\delta\theta}L^2}) & \longrightarrow (\mbox{\rm Ext}(S^\pm_{\overline X,\epsilon}),
\parallel\cdot\parallel_{e^{\delta\theta}L^2})\\
&    \overset{{\mathcal J}^\pm}{\longrightarrow}
\left[ {{\mathcal H}}_{S_{\overline M}}(]-\epsilon,\epsilon[)\oplus
{{\mathcal H}}_{S_{\overline M}}(]-\epsilon,\epsilon[)\right]^\pm
\end{split}
\end{equation*}

is continuous and exact in the middle.

\item[(b)] $h^\pm_{\Gamma, \epsilon}=\dim_\Gamma(\mbox{\rm im}({\mathcal J}^\pm)).$
\end{enumerate}
\end{lem}

\rightline{$\blacksquare$}

Let ${\mathcal J}:={\mathcal J}^++{\mathcal J}^-.$ We now want to show
that im$({\mathcal J})={{\mathcal V}}\oplus{\mathcal W}$ in 
${{\mathcal H}}_{S_{\overline M}}(]-\epsilon,\epsilon[)\oplus
{{\mathcal H}}_{S_{\overline M}}(]-\epsilon,\epsilon[).$ As $\tau_{\overline X}$
operates on im$({\mathcal J})$ and maps ${{\mathcal V}}$ and ${\mathcal W}$
to each other, the $\pm1$-Eigenspaces im$({\mathcal J}^\pm)$ must be isomorphic.
With a view of Lemma \ref{l5.8} and \ref{l5.18} this then implies the above Proposition.

To describe the structure of the image of ${\mathcal J}$, we consider the complement
of $L^2$-$\ker(S_{\overline X,\epsilon})$ in Ext$(S_{\overline X,\epsilon}).$ 
Denoting by $\perp$ the orhto-complement with regard to the nondegenerate pairing
 $e^{-\delta\theta}L^2\times e^{\delta\theta}L^2\to{\C},$ 
we see that

$$(L^2\mbox{\rm-ker}(S_{\overline X,\epsilon}),\parallel\cdot\parallel_
{e^{-\delta\theta}})^\perp=(e^{-\delta\theta}L^2\mbox{\rm-ker} (S_{\overline X,\epsilon}
))^\perp=\mbox{\rm cl}(e^{\delta\theta}L^2\mbox{\rm-im}(S_{\overline X,\epsilon})).$$

Writing $\K:=\mbox{\rm cl}(e^{\delta\theta}L^2\mbox{\rm-im}(S_{\overline X,\epsilon}))
\cap$ Ext$(S_{\overline X,\epsilon})$, it follows that im$({\mathcal J})=
{\mathcal J}(\K).$ The following Lemma is the reason why we have to look at the modification 
$S_{\overline X,\epsilon}$ instead of  $S_{\overline X}$ also in this Section.

\begin{lem}\label{l5.19} Let $0<\delta<\epsilon$.
\begin{enumerate}
\item[(a)] $e^{\delta\theta}L^2\mbox{\rm-im}(S_{\overline X,\epsilon})$ is $\Gamma$-dense
in $\mbox{\rm cl}(e^{\delta\theta}L^2\mbox{\rm-im}(S_{\overline X,\epsilon}))$, i.e. for every
$\kappa>0$ there is a closed subspace $\M\subset e^{\delta\theta}
L^2\mbox{\rm-im}(S_{\overline X,\epsilon})$ such that $\dim_\Gamma(\mbox{\rm cl}(e^{\delta\theta}
L^2\mbox{\rm-im}(S_{\overline X,\epsilon})))-\dim_\Gamma(\M)<\kappa.$

\item[(b)] $\K^0:=e^{\delta\theta}L^2\mbox{\rm-im}(S_{\overline X,\epsilon})\cap
\mbox{\rm Ext}(S_{\overline X,\epsilon})$ is $e^{\delta\theta}L^2$-dense in $\K.$
\end{enumerate}
\end{lem}

{\it Proof.} (a) is a direct consequence of the  $\Gamma$-Fredholm-property (Cf. (\ref{25}) and Section 6.2)
of the operator $T=S_{\overline X,\epsilon}$
on $e^{\delta\theta}L^2(\overline E)$. It implies that
$0$ is not in the  $\Gamma$-essential part of the spectrum of $TT^*$
(where $T^*$ denotes the adjoint of $T$ on $e^{\delta\theta}L^2$) and
the spaces 
$\M_\kappa:= \H_{TT^*}({\R}-]-\kappa,\kappa[)$ are of finite $\Gamma$-codimension
$\dim_\Gamma(\H_{TT^*}(]-\kappa,\kappa[-\{0\}))$ in $L^2\mbox{\rm-im}(T).$

(b) is a simple consequence of (a), cf. \cite{Sh2}.

\rightline{$\blacksquare$}

To show that ${\mathcal J}(\K)$ can be decomposed into a direct sum as claimed,
this Lemma and the continuity of the map $\mathcal J$ allow
to restrict our considerations to the dense subspace ${\mathcal J}(\K^0)$.

Thus, let $\xi\in\K^0.$ According to the definition of $\K^0$ there exists 
$\alpha\in e^{\delta\theta}L^2(\Lambda T^*\overline X)$ with $\xi=
S_{\overline X,\epsilon}\alpha$ and $S^2_{\overline X,\epsilon}\alpha=0.$
On the cylinder we can then write $\alpha=\alpha_0+\alpha_1 \wedge dr$ with
$\alpha_0,\alpha_1\in e^{\delta\theta}H^\infty(\overline Z,\Lambda T^*
\overline M)).$

Using the spectral resolution $V:L^2(\Lambda T^*\overline M)\to L^2({\R}_\lambda\times {\N},\mu_{S_{
\overline M}})$ of  $S_{\overline M}$, the $\alpha_l$ satisfy the following equation for $r>3$

$$-(\mbox{$\frac{\partial}{\partial r}$})^2 V\alpha_l+(1-\chi_{]\epsilon,
\epsilon[}(\lambda))\lambda^2 V\alpha_l=0.$$

The solution is of the general form

$$\alpha_l(x,r)=r\beta_{l,1}(x)r+\beta_{l,2}(x)+O(e^{-\epsilon r}),$$

for suitable $\beta_{l,1},\beta_{l,2}\in \H_{S_{\overline M}}(]-\epsilon,
\epsilon[).$ We can now write

$$ S_{\overline X,\epsilon}=d_{\overline X,
\epsilon}+d^*_{\overline X, \epsilon} \,\mbox{\rm with}\; 
d_{\overline X,\epsilon}=d_{\overline X}-d_{\overline M}\vartheta\widehat
\Pi_\epsilon\;\mbox{\rm iund}\; d^*_{\overline X,\epsilon}=d^*_{\overline X}-
d^*_{\overline M}\vartheta\widehat\Pi_\epsilon.$$

Then, in the area $r>3$ 

\begin{equation*}
\begin{split}
d_{\overline X,\epsilon}\alpha_0(x,r) & = (\epsilon(dr)\mbox{$\frac{\partial}
{\partial r}$}+d_{\overline M}(1-\widehat\Pi_\epsilon))(r\beta_{0,1}(x)+\beta_{0,
2}(x)+O(e^{-\epsilon r}))  \\
& = dr \wedge \beta_{0,1}(x)+O(e^{-\epsilon r}),
\end{split}
\end{equation*}

because $(1-\widehat\Pi_\epsilon)\beta_{0,j}=0.$ In the same manner one shows
$d_{\overline X,\epsilon}\alpha_1\wedge dr(x,r)=O(e^{-\epsilon r}).$
 The analogous calculation for $d^*_{\overline X,
\epsilon}\alpha$ gives

\begin{equation*}
\begin{split}
d^*_{\overline X,\epsilon}&\alpha_1(x,r)\wedge dr   \\
&= (-\iota(\mbox{$\frac{\partial}
{\partial r}$}) \mbox{$\frac{\partial}{\partial r}$}+d^*_{\overline M}
(1-\widehat\Pi_\epsilon))(r\beta_{1,1}(x)\wedge dr   +\beta_{1,2}(x)+O(e^{-\epsilon r})) \\
& = -(-1)^{|\beta_{1,1}|}\beta_{1,1}(x)+O(e^{-\epsilon r}),
\end{split}
\end{equation*}

and $d_{\overline X,\epsilon}\alpha_0(x,r)=O(e^{-\epsilon r})$.
We have shown
$${\mathcal J}(\xi)= {\mathcal J}(d_{\overline X,\epsilon}\alpha+d^*_{\overline X,\epsilon}\alpha)
=0\oplus (-1)^{|\beta_{0,1}|}\beta_{0,1}-(-1)^{|\beta_{1,1}|}
\beta_{1,1}\oplus 0,$$

that is, the image of $\mathcal J$ can be decomposed into a direct sum as claimed. This
finishes the proof of Proposition \ref{s5.15}.

\rightline{$\blacksquare$}

\appendix

\section{Clifford algebra conventions}

Denote by ${\C}(k)$ the (complex) Clifford algebra over the euclidean space
${\R}^k,$ with generators ${\bold c}_1,\ldots,{\bold c}_k$ satisfying 
${\bold c}_i{\bold c}_j+{\bold c}_j{\bold c}_i=-2\delta_{ij}.$ The algebra ${\C}(k)$
is ${\Z}_2$-graded: ${\C}(k)={\C}^+(k)\oplus 
{\C}^-(k),$ and the map ${\bold c}_i\longmapsto{\bold c}_i
{\bold c}_{k+1}$ defines an isomorphism ${\C}(k)\overset{\sim}{\to}
{\C}^+(k+1)$. The  volume element $\tau_k:=i^{[(k+1)/2]}{\bold c}_1
\ldots{\bold c}_k\in{\C}(k)$ satisfies $\tau^2_k=1$ and thus induces a ${\Z}_2$-grading 
on the representations of ${\C}(k).$
Note however, that $\tau_k{\bold c}=-(-1)^k{\bold c}\tau_k$ for ${\bold c}\in{\R}^k\subset 
{\C}(k)$ implies that this grading is trivial on  irreducible representations ${\C}(k)$, when $k$ is odd.

${\C}(2l)$ has a unique irreducible representation, called its spinor space and denoted by
$S(2l).$ Its dimension is $\dim S(2l)=2^l.$ Decomposing into the  $\pm1$-Eigenspaces
of $\tau_{2l}$ we write $S(2l)=S^+(2l)\oplus S^-(2l).$ Via the identification ${\C}
(2l-1)\cong {\C}^+(2l)$ the spaces $S^+(2l)$, $S^-(2l)$ are non-equivalent
irreducible representations of ${\C}(2l-1)-$, which can be considered as being
isomorphic representations of ${\C}(2l-2)\cong {\C}^+(2l-1)$ via the map $S^+(2l)
\overset{c_{2l}}{\to}S^-(2l)$. This of course is then just the representation $S(2l-2)$
of  ${\C}(2l-2)$.
For  $S^\pm(2l)$ we also write  $S^\pm(2l-1)$ when these spaces are seen as 
representations of ${\C}(2l-1)$.

It is easy to see that ${\C}(2l)$ acts injectively on $S(2l)$. 
Comparison of dimensions then yields ${\C}(2l)\cong End(S(2l)),$ and,
using ${\C}(2l-1)\cong{\C}^+(2l)$ also ${\C}(2l-1)\cong 
End^+(S(2l)).$ The identification ${\C}(2l-1)\to End(S^\pm(2l-1))$ maps
$\tau_{2l-1}$ to $\pm 1$ and thus has null space $(1\mp\tau_{2l-1}){\C}(2l-1).$

The traces $tr^\pm$ on $End(S^\pm(2l-1))$ and the graded trace $str$ on
$End(S(2l))$ then induce traces on ${\C}(2l-1)$ and 
${\C}(2l).$  On elements of the form ${\bold c}_I:={\bold c}_{i1}
\ldots{\bold c}_{i|I|}$ where $I=\{i_1\le\ldots\le i_{|I|}\}\subset \{1,\ldots,
k\}$ these are calculated as follows

\begin{lem}\label{A.1} 
\begin{enumerate}
\item[(a)] In ${\C}(2l)$ we have $\mbox{\rm str}(\tau_{2l})=2^l$ 
and $\mbox{\rm str}(1)=\mbox{\rm str}({\bold c}_I)=0$ for $I\neq\{1,\ldots,k\}.$

\item[(b)] In ${\C}(2l-1)$ we have $str^+(\tau_{2l-1})=-tr^-(\tau_{2l-1})=
tr^\pm(1)=2^{l-1}$ and for $I\neq\{1,\ldots,k\}$ we have $tr^\pm({\bold c}_1)=
0.$
\end{enumerate}
\par

On $({\C}(2l-1)-{\C})\subset{\C}(2l)$ therefore $tr^\pm
(\bullet
)=\mp\frac{1}{2} \mbox{\rm str}({\bold c}_{2l}\bullet)$ and on $({\C}(2l)\subset
{\C}(2l+1)$ we have $\mbox{\rm str}(\bullet)=\pm i tr^\pm({\bold c}_{2l+1}\bullet)$
\end{lem}

{\it Proof.} Cf. \cite{BGV}, Proposition 3.21

\rightline{$\blacksquare$}

The map $S^+(2l)\overset{c_{2l}}{\to} S^-(2l)$ gives an identification
$S(2l)\cong S^\pm(2l-1)\oplus S^\pm(2l-1).$ In this representation,  ${\C}(2l)$ 
acts on $S(2l)$ as follows

\begin{equation*}
\begin{split}
&  {\bold c}_i\in{\C}(2l-1)\overset{\wedge}{=}
\begin{pmatrix}
0  &  \pm{\bold c}_i  \\
\pm {\bold c}_i  &  0
\end{pmatrix}\;
{\bold c}_{2l}\overset{\wedge}{=}
\begin{pmatrix}
0  &  -1  \\
1  &  0
\end{pmatrix}\\
&  \mbox{and}\;
str
\begin{pmatrix}
\phi_1 & \phi_2 \\
\phi_3 & \phi_4
\end{pmatrix}=tr^\pm(\phi_1)-tr^\pm(\phi_4)
\end{split}
\end{equation*}

\end{document}